%% file: main.tex
\documentclass[12pt,twoside,english]{article}
\usepackage{ifthen}                % if then else Abfragen
\pdfoutput=1
\newboolean{boldeqnitalic}
\setboolean{boldeqnitalic}{true}

\usepackage{caption}
\usepackage{umlaute}
\usepackage[intlimits] {amsmath}   % muss for defjs# stehen 
\usepackage{defjs2}                % eigene Definitionen    
\usepackage{epsfig}
\usepackage{psfig}
\usepackage{graphicx}
\usepackage{psfrag}
\usepackage[usenames]{color}
\usepackage{colortbl}
\usepackage{ifthen}
\usepackage{fancyhdr}
\usepackage{rotating}
\usepackage{multirow}
\usepackage{booktabs}              % Dicke Tabellenlinien
\usepackage{amssymb}
\usepackage{amsbsy}
\usepackage{bm}                    % fette Schrift in Formeln
\usepackage{babel}
\usepackage{theorem}
\usepackage{upgreek}  % gerade (roman) griechische Buchstaben z.B. \upalpha
\usepackage{pifont}   % zusaetzliche Schriften, z.B. eingekreiste Zahlen
\usepackage[ruled]{algorithm2e}
\usepackage{algpseudocode}
\usepackage{mathrsfs}
\usepackage{datetime}
\usepackage[all]{xy}
\usepackage{rotating}
\usepackage{tabulary}
\usepackage{enumitem}
\usepackage{tikz}

\definecolor{dunkelgrau}{rgb}{0.8,0.8,0.8}
\definecolor{hellgrau}{rgb}{0.90,0.90,0.90} %... slightly darker 
\fancypagestyle{plain}{%
  \fancyhead{}%
  \fancyfoot[c]{\sffamily\thepage}%
}
\makeatletter                           % Leere Fuellseiten
\def\cleardoublepage{\clearpage\if@twoside \ifodd\c@page\else
  \hbox{}
  \vspace*{\fill}
  \thispagestyle{empty}
  \newpage
  \if@twocolumn\hbox{}\newpage\fi\fi\fi}
\makeatother
\begin{document}
\unitlength1.0cm
\frenchspacing
\include{antrag_text}

% ------- layout-datei --------------
% Variante I - elsevier, numbers, references with URLs
%\bibliographystyle{elsarticle-num}

% Variante II - numbers, references without URLs
\bibliographystyle{plainnat}

%\bibliographystyle{plaindin}
%\bibliographystyle{plaindin_shortname2}
%\bibliographystyle{elsarticle-num-names}
% ------- bib-datei --------------

%\begin{appendix}
%\include{appendixA}
%\end{appendix}
\end{document}

%% file: antrag_text.tex
\thispagestyle{empty}

\input{Frontpage-1} 

\bigskip
 
\section{Introduction}
\label{sec:intro}
For the solution of homogenization problems, two-scale finite element methods like the FE$^2$ method and the Finite Element Heterogeneous Multiscale Method FE-HMM combine in a favorable manner the needs for both accuracy and efficiency in that they sample the real microstructure in representative volume elements RVEs of confined size. By virtue of the equivalence of the coarse scale energy density with the fine-scale energy density these methods stand on sound ground in physics, by virtue of the a priori estimates developed for FE-HMM, they are equipped with necessary ingredients in mathematics. 

In spite of considerable achievements in the field, there is space for improvements; the present work addresses three of them. 
\begin{enumerate}
	\item {\bf Nonlinear FE-HMM}. A nonlinear version of the Finite Element Heterogeneous Multiscale Method FE-HMM is developed for the homogenization of solids. 
	\\[-7mm] 
	\item {\bf A priori estimates}. An assessment of the a priori estimates of linear FE-HMM is carried out in the nonlinear regime of hyperelasticity. 
	\\[-7mm]
	\item {\bf Solution algorithm}. The revision and improvement of the standard algorithmic structure of nonlinear two-scale finite element methods for homogenization like the FE$^2$ method and the present nonlinear FE-HMM. 
\end{enumerate}

These aims shall be motivated and put into perspective by a brief account of the state of the art. 

\subsection{Development of a nonlinear FE-HMM}
 
FE-HMM formulations for problems in solid mechanics are so far restricted to linear problems, i.e. to linear elasticity in a geometrical linear setting; \cite{Assyr2006} presents a fully discrete convergence analysis taking into account the discretization errors at both micro and macro levels. 

FE$^2$ and FE-HMM are equivalent methods, for a comparison see \cite{EidelFischer2018}; they are siblings (if not twins) having different parents. FE$^2$, an engineering mechanics off-spring, is quite outgoing since its infancy and, unburdened from the necessity of delivering mathematical results on existence and uniqueness of solutions, has appreciated direct contact with nonlinear problems since its first appearance  \cite{MichelMoulinecSuquet1999}, \cite{Miehe-etal-1999a}, \cite{Miehe-etal-1999b}, \cite{FeyelChaboche2000}, \cite{Kouznetsova-etal2001}. FE-HMM as a method of mathematicians has grown up subject to the restriction, first prove, then apply; consequently its playground of applications was mostly in linear problems, comprehensive overviews are given in \cite{E-Engquist-Li-Ren-VandenEijnden2007}, \cite{Assyr2009}, \cite{Assyr-etal2012}. For the mathematical analysis of FE-HMM in the nonlinear regime of elliptic PDEs, results are sparse. \cite{AbdulleHuber2016} consider nonlinear monotone elliptic problems. For the same case \cite{HenningOhlberger2015} derive a posteriori error estimates for the $L^2$-error between the HMM approximation and the solution of the macroscopic limit equation. \cite{Nejad-Wieners-2019} present an FE-HMM formulation for inelasticity that draws on the tangent moduli for stiffness-transfer and therein very much resembles the FE$^2$ method.

The present work proposes a nonlinear FE-HMM formulation for solid mechanics, which implies geometrical nonlinearity, either along with linear elasticity or nonlinear hyperelasticity. The extension to inelastic constitutive laws is a minor step following the perspective that inelasticity --in finite element numerics-- can be understood as a type of nonlinear elasticity along with evolution equations for the inelastic deformation part treated as internal variables. 
  
\subsection{Test of the a priori estimates of linear FE-HMM for hyperelasticity}
Error estimates are a necessary equipment of numerical methods to serve as reliable tools in science and engineering. By virtue of its foundation in mathematical homogenization by asymptotic expansion a priori estimates are available for FE-HMM; the elliptic case is considered in \cite{E-Ming-Zhang-2005}, \cite{Ohlberger2005}, the elliptic case of linear elasticity in a geometrical linear setting in \cite{Assyr2006}, \cite{Assyr2009}. 
An assessment of the a priori estimates for various energetically consistent boundary conditions applied to the micro domain is presented in \cite{EidelFischer2018a}.

In the present work we assess the convergence of FE-HMM solutions in the nonlinear regime of hyperelastic material behavior, although the estimates are derived for the setting of linear elasticity (generally, for linear elliptic PDEs). This will be carried out by a comparison with the linear case, which clarifies, whether the deviation from the estimates are caused by a lowered regularity of the boundary value problem or whether the deviation follows from the nonlinearities of the considered setting. 
 
\subsection{Speedup through modified algorithmic structure of two-scale FEM}
In view of the considerable computational costs of two-scale FEM in homogenization there is a need for methodic advancements aiming at improved performance. For that aim several new and modified methods have been developed beyond parallelization.

Homogenization techniques based on Fast Fourier Transforms (FFT) were introduced in \cite{Moulinec-Suquet-1994}, \cite{Moulinec-Suquet-1998}, which are built on a Lippmann-Schwinger formulation \cite{Zeller-Dederichs}, \cite{Kroener-BOOK-1971}; applicability and efficiency was demonstrated in \cite{Schneider-Ospald-Kabel-2015}, \cite{Schneider-Merkert-Kabel-2017}, and many more.   
Reduced basis methods and model order reduction exhibit considerable potential to speed up simulations; \cite{Yvonet-He-2007} introduced a general (Galerkin) projection based reduced order model referred to as R3M where the reduced basis is obtained via snapshot proper orthogonal decomposition. Applications to nonlinear elasticity in \cite{Radermacher-Reese}, to homogenization problems in \cite{Soldner-etal2017}. The derived reduced problems are low-dimensional and nonlinear. 
In an effort to reduce the computational complexity by a modified modeling \cite{SchroederBalzaniBrands2010} consider in FE$^2$ statistically similar representative volume elements (SSRVE) in order to replace the true microstructure in its full geometrical complexity by a simplified surrogate that resembles the original one by geometrical features as analyzed by different geometrical similarity measures. Domain-decomposition based on a FETI (finite element tearing and interconnect) have been proposed and successfully applied \cite{RheinbachKlawonn2006}, \cite{RheinbachKlawonn2010}. 

Pixel- and voxel-based representations of reconstructed  microstructures obtained from tomographic imaging methods is frequently used in various fields \cite{Holzer.2011},  \cite{Saenger-etal-2011}, \cite{Andrae-etal-2013}. Suchlike spatially uniform discretizations are typically the point of departure for quadtree- and octree-based mesh coarsening which can strongly reduce computational efforts in homogenization \cite{Legrain-etal2011}, \cite{Lian-etal2013}. While the gain in efficiency was quantified earlier, the error introduced by adaptive mesh coarsening was for the first time assessed in \cite{FischerEidel2019}.  
 
The present work aims at speeding up computations through an improvement of the standard algorithmic concept of FE$^2$. The discretized variational form of the balance of linear momentum both on the macro level and the micro level amount to systems of nonlinear algebraic equations, which are each typically but not necessarily solved by Newton's method. More important, the algorithmic solution framework is not monolitic, but exhibits a staggered structure, which provides a mutual exchange of coupling quantities (micro-to-macro: current micro stiffness and stress, macro-to-micro: current macro deformation) between the macro problem and the micro problems. 

The very standard of suchlike staggered scheme is realized by a nested loop embedding the full solution of the micro problems into one macro solution iteration step as sketched in the pseudocode of ''Algorithm 1''. 

The attribute \emph{standard} is set by a multitude of contributions to the field, \cite{FeyelChaboche2000} p.~313, \cite{Kouznetsova-etal2001} pp.~41,42, \cite{Miehe-Koch-2002} pp.~302,303, \cite{Miehe-2003} p.~574, \cite{Feyel2003} p.~3236, \cite{Kouznetsova-etal2004} p.~5532, for thermo-mechanical problems \cite{OezdemirBrekelmansGeers2008} p.~606, \cite{TemizerWriggers2008} p.~499, \cite{Schroeder2014} p.~8, for different coupling conditions at finite deformations \cite{Saeb-Steinmann-Javili-2016} p.~20, for nonlinear monotone elliptic problems in FE-HMM \cite{AbdulleHuber2016} p.~966, \cite{KlawonnKohelerLanserRheinbach2019} p.~7, to name but a few.

\begin{minipage}[t]{16cm}
	\centering
	\begin{minipage}[t]{10.0cm}
		\begin{algorithm}[H]
			\While{macro residual $<$ macro tolerance}{
				macro iteration\;
				\While{micro residual $<$ micro tolerance}{
					micro iteration\;      
				}
			}
			\caption{nested (standard)}
		\end{algorithm}
	\end{minipage}
\end{minipage} 

In the present work a revision of this established scheme is carried out along with the proposal to replace it by an algorithm of direct alternations between micro and macro iterations as sketched in the pseudocode box ''Algorithm 2''. 

\begin{minipage}[t]{16cm}
	\centering
	\begin{minipage}[t]{9.5cm}
		\begin{algorithm}[H]
			\While{macro residual $<$ macro tolerance}{
				macro iteration\;
				micro iteration\;      
			}
			\caption{alternating (novel)}
		\end{algorithm} 
	\end{minipage}
\end{minipage}

The novel scheme, which was introduced in \cite{Eidel-etal-2018}, is described in more detail in Sec.~\ref{subsec:novel}. Its properties and benefits in terms of accuracy, efficiency and stability with respect to large load step sizes are assessed in the example section \ref{sec:NumericalExamples}.  
 
\bigskip

The paper is organized along the route of the above aims.   

\section{Strong and variational forms}

We consider a body $\mathcal{B}_0$, a bounded subset of $\mathbb{R}^{n_{\text{dim}}}$, ${n_{\text{dim}}}= 2,3$, with boundary 
$\partial \mathcal{B}_0= \partial \mathcal{B}_{0D} \cup \partial \mathcal{B}_{0N}$ where $\partial\mathcal{B}_{0D}$ and $\partial\mathcal{B}_{0N}$ are disjoint sets.
The closure of the body $\mathcal{B}$ is denoted by $\overline{\mathcal{B}}$.
The body undergoes deformation $\bm \varphi:\Omega \rightarrow \mathbb{R}^{n_{\text{dim}}}$ with deformation gradient $\bm F = \partial_{\bm X} \bm \varphi(\bm X)$ and the Jacobian $J=\text{det}\bm F>0$, where $\bm X$ is a material point in the reference configuration. The body is subject to body forces $\bm b$ and surface tractions $\bm g$ and in static equilibrium 

\begin{equation}
\text{Div}[\bm P] + \rho \, \bm b = \bm 0
\label{eq:Balance-linear-momentum}
\end{equation}
in terms of the first Piola-Kirchhoff stress tensor $\bm P$, density $\rho$ in the reference configuration, and body forces $\bm b$ neglecting inertia terms.

Dirichlet and Neumann boundary conditions are prescribed by
\begin{equation}
\label{eq:BCs}
\bm u = \hat{\bm u} \quad \text{on} \quad \partial\mathcal{B}_{0D} \qquad \text{and} \quad
\bm P \cdot \bm N = \bm t \quad \text{on} \quad \partial\mathcal{B}_{0N} \, .
\end{equation}

The corresponding variational form of the Boundary Value Problem (BVP) \eqref{eq:Balance-linear-momentum}, \eqref{eq:BCs} reads 

\begin{equation}
G:= \int_{\mathcal{B}_0} (\text{Div}[\bm P] + \rho \, \bm b) \cdot \delta \bm u \, \text{d}V = 0\, ,
\end{equation}
with virtual displacements/test functions $\delta \bm u$ which can be transformed by $\text{Div}[\bm P^T \cdot \delta \bm u] - \bm P : \text{Grad}[\bm P^T \cdot \delta \bm u]$
into 
 
\begin{equation}
G:= \underbrace{\int_{\mathcal{B}_0} \bm P : \text{Grad}[\delta \bm u] \, \text{d}V}_{\displaystyle = G^{\text{ext}}}
    \underbrace{- \int_{\partial \mathcal{B}_{0N}} \bm t \cdot \delta \bm u \, \text{d}A
     - \int_{\mathcal{B}_0} \rho \, \bm b \cdot \delta \bm u \, \text{d}V}_{\displaystyle = G^{\text{int}}} = 0\, ,
\end{equation}  

which has to hold for all $\delta \bm u \in \mathcal{V}$, where $\mathcal{V}$ is the space of admissible displacements, i.e. virtual displacements that fulfill homogeneous Dirichlet boundary conditions
\begin{equation}
\label{eq:HilbertSpaceV}
\mathcal{V}=\{\delta \bm u; \bm u \in H^1(\mathcal{B}_0)^{{n_{\text{dim}}}}, \bm u|_{\partial \mathcal{B}_{0D}} = \bm 0 \} \, .
\end{equation}

With  
\begin{equation}
\bm P : \text{Grad}[\delta \bm u] = \bm S : \bm F^T \cdot \text{Grad}[\delta \bm u]
= \bm S : \dfrac{1}{2} (\bm F^T \cdot \delta \bm F +\delta \bm F^T \cdot \bm F)
= \bm S : \delta \bm E
\end{equation}

the virtual work can be equally expressed by the second Piola-Kirchhoff stress tensor $\bm S$ and the Green-Lagrange strain tensor $\bm E$
\begin{equation}
\label{eq:Variational-form-Lagrange-continuum}
G:= \int_{\mathcal{B}_0} \dfrac{1}{2}\bm S : \delta \bm E \, \text{d}V 
   - \int_{\partial \mathcal{B}_{0N}} \bm t \cdot \delta \bm u \, \text{d}A
	- \int_{\mathcal{B}_0} \rho \, \bm b \cdot \delta \bm u \, \text{d}V  = 0\, .
\end{equation}   

\subsection{Characteristics and conditions for computational homogenization}
\label{subsec:Charactistics-and-conditions-for-CompHomog}
For the solution of the BVP in its variational form, a constitutive law  must be defined in the domain $\mathcal{B}_0$ for which we assume matter to be heterogeneous in its properties, such that the solution to the BVP shall be found through homogenization. For the application of homogenization in terms of computational micro-macro transitions several conditions must be fulfilled, which partly refer to the characteristics of the physical problem, partly to mathematical homogenization, and partly to the properties of the solution method. Consequently, the following assumptions are made:
\begin{itemize}
	\item {\bf Scale separation}. A prerequisite of homogenization is that the characteristic length of microstructural features must be \emph{considerably} smaller than the characteristic length on the macro scale\footnote{We use macro/coarse scale and micro/fine scale each as synonyms.} $L_{\text{mac}} \gg L_{\text{mic}}$. The deeper root of this requirement is the separation of variables in mathematical homogenization which requires a separation of length scales with the implicit, tacit or explicit assumption that the micro solution and macro solution equally reflect this scale separation \cite{Bensoussan-Lions-Papanicolau-BOOK-1976}, \cite{Sanchez-Palencia-BOOK-1980},  \cite{GuedesKikuchi1990}, \cite{Allaire1992}, \cite{Cioranescu-Donato-BOOK-1999}.
	\item {\bf Sampling domains}. For the microstructure on the fine scale a statistically representative volume element RVE is available, \cite{Hill1963}, \cite{Ostoja-Starzewski2006}, \cite{Sab-1992}, \cite{Drugan-Willis-1996}.  
	\item {\bf Constitutive equations}. The type of constitutive law and its parameters are available for matter in the RVE, for its constituents (grains, phases, etc.) and  interfaces, but not on the macro scale in terms of effective laws and their parameters.  
	\item {\bf Coupling conditions}. The coupling conditions cover a kinematical and an energetic aspect. The equivalence of energy densities of the coarse scale with the fine scale is a cornerstone of homogenization, from which different boundary conditions on the RVE are derived each fulfilling this equivalence \cite{Hill1963}, \cite{Mandel-BOOK-1971}, \cite{Hill1972}, \cite{E-Engquist-2003}. The kinematical embedding of the RVE into the macro domain is part of the modeling assumptions; if the deformation applied to the RVE is homogeneous, it is referred to as first order computational homogenization.    
	\item {\bf Scale interaction and data transfer}. The conditions for these two items follow from the above items of constitutive equations and coupling conditions. The macro BVP is the driver of the micro BVP through deformation in the coupling conditions, the top-down direction of data transfer. Since constitutive equations with parameters are a priori exclusively available on the micro scale, the numerical solution of the macro BVP requires stiffness and stress from the fine scale, the bottom-up direction of data transfer.  
\end{itemize} 
\begin{Figure}[htbp]
	\begin{minipage}{16.5cm}  
		\hspace*{16mm}
		\includegraphics[width=12.0cm, angle=0, clip=]{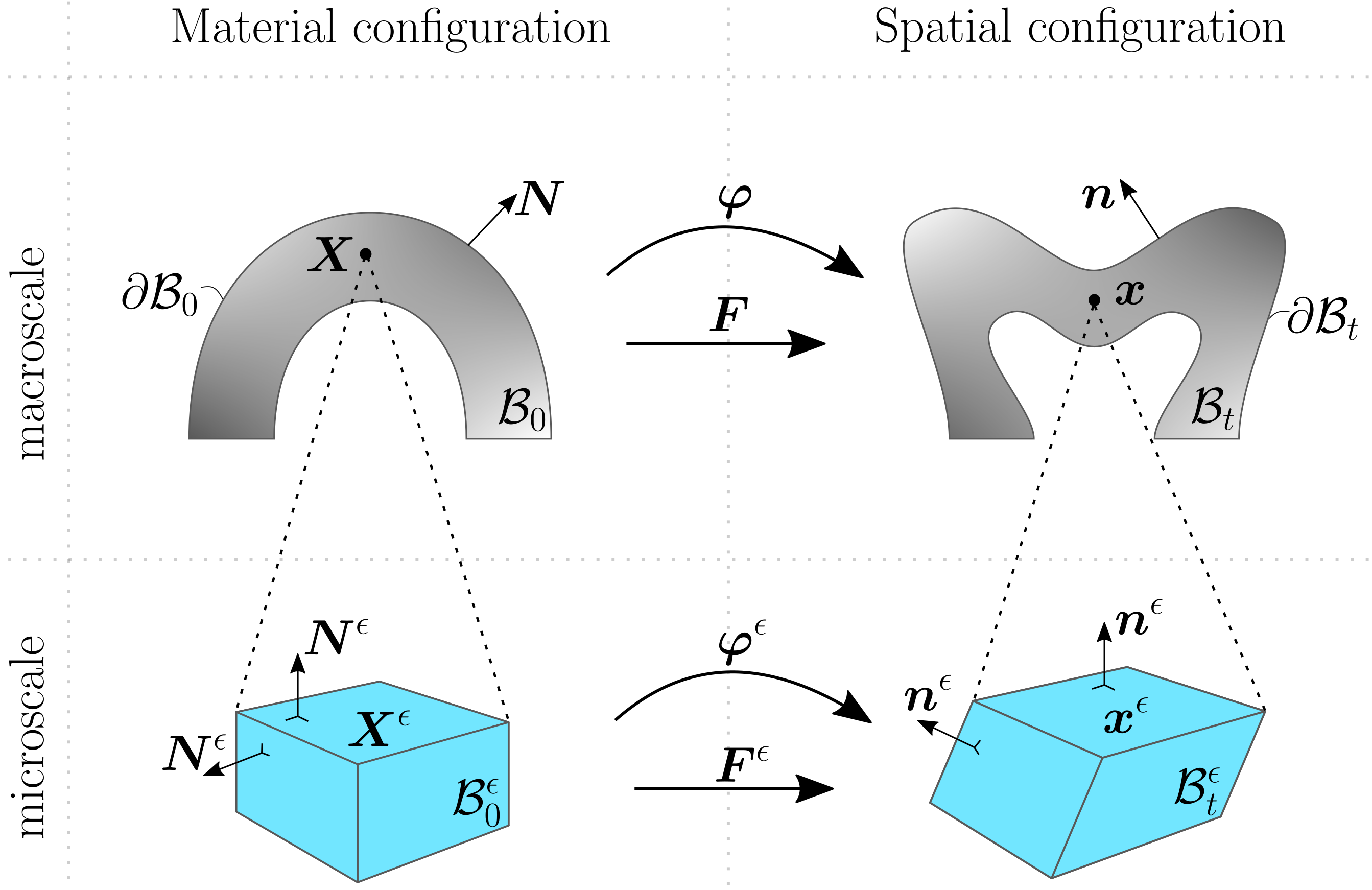}  
	\end{minipage}
	\caption{{\bf Homogenization}: The general setting of micro-macro transition in material and spatial configurations. Sketch adapted from \cite{Saeb-Steinmann-Javili-2016}.} 
	\label{fig:Homogenization-scheme} 
\end{Figure}

While balance laws are scale invariant and therefore are known and apply for the macro as well as the micro scale, the effective constitutive law along with its parameters is not known for the macro scale but exclusively for the micro scale, here in terms of the hyperelasticity relations  
\begin{equation}
\bm{P}^{\epsilon} \ = \ \dfrac{\partial \psi^{\epsilon}}{\partial \bm F^{\epsilon}} \, , \quad 
\mathbb{A}^{\epsilon} \ = \ \dfrac{\partial^2 \psi^{\epsilon}}{\partial \bm F^{\epsilon}\partial \bm F^{\epsilon}} \, , 
\quad 
\bm S^{\epsilon} = \bm F^{\epsilon \, -1} \, \bm P^{\epsilon} \, , 
\label{eq:hyperelasticity-PF}
\end{equation} 
where we indicate micro scale quantities by $\epsilon$ as borrowed from mathematical homogenization, where it refers to the (periodic) cell size.   

The Hill-Mandel energy equivalence principle is formulated as an equivalence of stress power in terms of the first Piola-Kirchhoff stress tensor and the material time derivative of the work-conjugate deformation gradient $\bm F$
\begin{equation}
\label{eq:Hill-Mandel}
  \bm P : \dot{\bm F} = \dfrac{1}{|\mathcal{B}_0|} \int_{\mathcal{B}_0} \bm P^{\epsilon} : \dot{\bm F}^{\epsilon} \, \text{d}V \, ,
\end{equation}
where their averages are defined according to 
\begin{equation}
  \label{eq:Average-of-P-and-F}
    \bm P \ = \ \langle \bm P^{\epsilon} \rangle_0 \, ,
    \qquad 
    \bm F \ = \ \langle \bm F^{\epsilon} \rangle_0 \, ,
    \quad \text{with} \quad
    \langle \{\bullet\} \rangle_0 = \dfrac{1}{|\mathcal{B}_0|}\int_{\mathcal{B}_0} \{\bullet\} \, \text{d}V \, ,
\end{equation}
with $|\mathcal{B}_0|$ the measure of $\mathcal{B}_0$. The volume average for the deformation gradient according to \eqref{eq:Average-of-P-and-F}$_2$ is only valid for microstructures without hollow spaces; for materials with voids, cracks, or for cellular materials it holds
\begin{equation}
 \label{eq:Average-of-F-precisely}
  \bm F := \dfrac{1}{|\mathcal{B}_0|} \int_{\partial \mathcal{B}_0} \bm x \otimes \bm N \, \text{d}A = \dfrac{1}{|\mathcal{B}_0|} \left[ \int_{\mathcal{B}_0} \bm F \, \text{d}V 
    - \int_{|\circ|} \bm x \otimes \bm N \, \text{d}A \right]\, ,
\end{equation}
where $\bm x$ is the position vector of a point in $\mathcal{B}_0$, $\bm N$ is the normal on $\partial {\mathcal B}_0$, and where $|\circ|$ represents the boundaries of hollow space.

Since the present FE-HMM framework is formulated in the reference configuration with the work conjugate pair of the second Piola-Kirchhoff stress and the Green-Lagrange strain tensor, the macro average stress must be consistent with \eqref{eq:Average-of-P-and-F} in \eqref{eq:Hill-Mandel}, it holds 
\begin{equation} 
{\color{black} \bm S \ = \ \langle \bm F^{\epsilon} \rangle^{-1}_0 \langle \bm P^{\epsilon} \rangle_0} \, .
\label{eq:macro-mean-stress-S} 
\end{equation} 
An averaged macro stress according to $\bm S^{\ast} :=\langle \bm S^{\epsilon} \rangle_0$ is not consistent with \eqref{eq:Average-of-P-and-F} since $\bm S \neq \bm S^{\ast}$.  

A framework for equivalence relationships for finite strain multi-scale solid constitutive models based on the volume averaging of the microscopic stress and deformation gradient fields over a representative volume element (RVE) is presented in \cite{SouzaNeto-Feijoo2008}. Based on a purely kinematically-based variational framework they derive sufficient conditions under which the volume average of the microscopic first Piola-Kirchhoff stress over the material configuration of the RVE is mechanically equivalent to the average of the microscopic Cauchy stress field over the spatial  configuration.  
  
Similar to the present work, the weak form in $\bm S$ and $\bm E$ in the context of FE$^2$ is chosen in \cite{GrytzMeschke2008} and \cite{vDijk-2016}.
 
%--------------------------------------------------------------------------------------------------------
\section{The Heterogeneous Multiscale Finite Element Method}
%--------------------------------------------------------------------------------------------------------

A nonlinear version of FE-HMM is derived for the application to geometrical nonlinear problems in solid mechanics based on linear elastic or (nonlinear) hyperelastic constitutive laws. In the example section \ref{sec:NumericalExamples} we will consider Neo-Hookean-type hyperelasticity. While the majority of homogenization schemes is formulated in reference configuration in terms of $\bm P$ and $\bm F$, we employ the work-conjugate pair $\bm S$ and $\bm E$. 

\subsection{Discretization and linearization of the variational form}
\label{subsec:Variational-FE-HMM-macro}

\begin{Figure}[htbp]
	\begin{minipage}{16.5cm}  
	 \hspace*{2mm}
     	\includegraphics[width=15.0cm, angle=0, clip=]{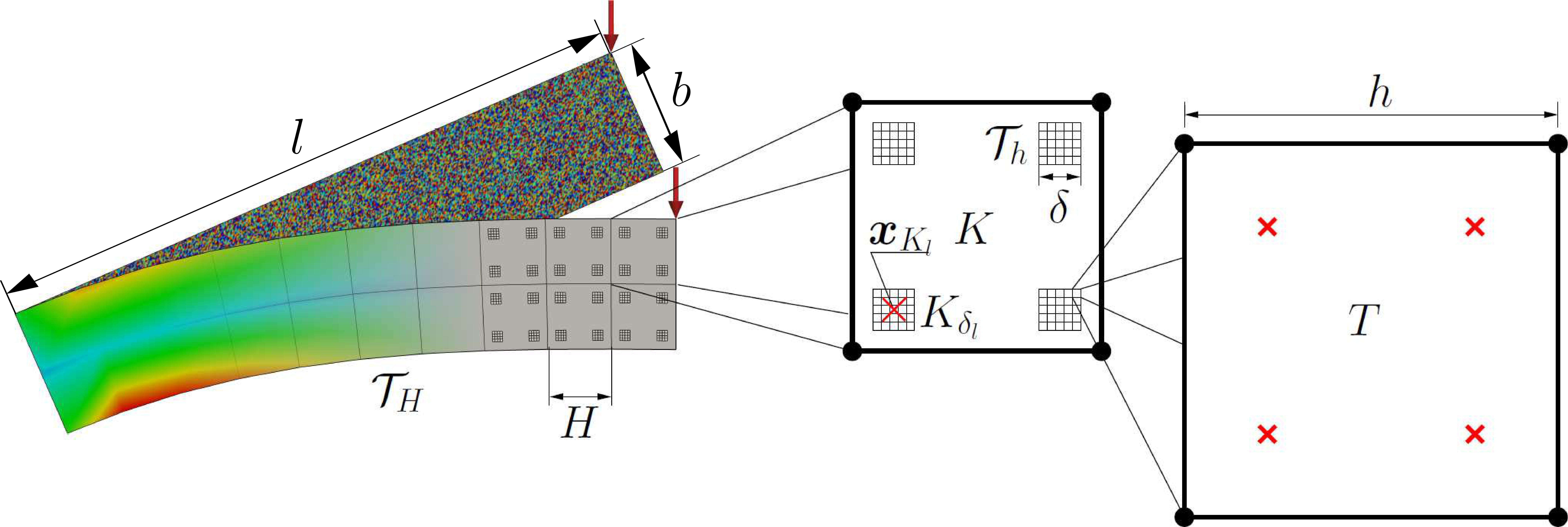}   
	\end{minipage}
	\caption{{\bf Micro-macro transition of FE-HMM}: The triangulation of the macro problem ${\cal T}_{H}$ results in finite elements $K$ of size $H$. Micro sampling domains of characteristic length $\delta$ in triangulation ${\cal T}_{h}$ into elements $T$ of size $h$ are attached to the macro quadrature points $\bm x_{K_{\delta_l}}$ of $K$, hence spanning the sampling domain $K_{\delta_l}=\bm x_{K_{\delta_l}} + \delta \, [-1/2, +1/2]^{n_{\text{dim}}}$, $\delta \geq \epsilon$, centered at quadrature points ($qp$), $l=1, \ldots, N_{qp}$. The separation of length scales imposes the condition $\{l,b\}\gg \epsilon$.} 
	\label{fig:MicMac-Problem-Meshing-etc} 
\end{Figure}
 
Here we consider the piecewise linear continuous FEM in macro and {\color{black} micro space}, respectively. 

We define a macro finite element space as 
\begin{equation}
 \mathcal{S}^p_{\partial \mathcal{B}_{0D}}(\mathcal{B}_0, {\mathcal T}_H) = \left\{ \bm u^H \in H^1(\mathcal{B}_0)^{n_{\text{dim}}}; \bm u^H|_{\partial \mathcal{B}_{0D}} = \bm 0; \bm u^H|_{K} \in {\mathcal{P}}^{p}(K)^{n_{\text{dim}}}, \, \forall \, K \in {\cal T}_{H} \right\}
 \label{eq:MacroFESpace}
\end{equation}
where ${\mathcal P}^{p}$ is the space of (in the present work: linear, $p=1$) polynomials on the element $K$, ${\mathcal T}_H$ the (quasi-uniform) triangulation of $\mathcal{B}_0\, \subset \, \mathbb{R}^{n_{\text{dim}}}$. The space $\mathcal{S}^{p}_{\partial \mathcal{B}_{0D}}$ is a subspace of $\mathcal{V}$ defined in \eqref{eq:HilbertSpaceV}. The characteristic element size on the macro scale $H$ (Fig.~\ref{fig:MicMac-Problem-Meshing-etc}) refers in index or superscript to macro scale quantities in the discretized setting, the characteristic micro element size $h$ analogously marks micro scale quantities of the discrete case.

For the solution of the macro scale BVP we use the two-scale FEM framework of the FE-HMM, as originally proposed in \cite{E-Engquist-2003} and analyzed for elliptic PDEs in \cite{E-Ming-Zhang-2005}.

The macro solution of the FE-HMM is given by the following variational form:

Find $\bm u^H \in \mathcal{S}_{\mathcal{B}_{0D}}(\mathcal{B}_0, \mathcal{T}_H)$ such that \begin{equation}
\label{eq:VariationalFormulationHMM}
 {\color{black}G^H(\bm u^H, \delta \bm u^H) = \int_{\mathcal{B}_0} \dfrac{1}{2}\bm S^H :\delta \bm E^H \,\text{d}V} 
            - \int_{\partial \mathcal{B}_{0N}} \bm t \ \cdot \ \delta \bm u^H \, \text{d}A
            - \int_{\mathcal{B}_0} \bm b \ \cdot \ \delta \bm u^H  \, \text{d}V
              \quad \forall \ \delta \bm u^H \in \mathcal S_{\partial \mathcal{B}_{0D}} (\mathcal{B}_0, \mathcal{T}_H) \, .
\end{equation}  
Standard finite element shape functions $N_I$ are used for the interpolation 
\begin{equation}
\label{eq:u^H,u^h}
\bm{u}^H = \ \sum_{I=1}^{N_{\text{node}}} N_I \bm{d}_I^H \, , \quad  
\bm{u}^h = \ \sum_{I=1}^{n_{\text{node}}} N_I \bm{d}_I^h \, ,
\end{equation}
of nodal displacement vectors $\bm{d}_I^H$, $\bm{d}_I^h$, where $N/n_{\text{node}}$ is the number of nodes per element on the macro/micro scale. Virtual displacements $\delta \bm{u}^H$ and $\delta \bm{u}^h$ are equally obtained through interpolation of nodal values $\delta \bm{d}_I^H$ and $\delta \bm{d}_I^h$, respectively. Consequently, we obtain for $G^H$ in   \eqref{eq:VariationalFormulationHMM}
 \begin{equation}
 G^H(\bm{u}^H, \delta \bm{u}^H) \ = \ \sum_{I=1}^{N_{\text{node}}}  \left(\delta \bm{d}_I^{H}\right)^T \bm{f}_I^H
 \quad \text{with} \quad 
  \bm{f}_I^H = \bm{f}_I^{\text{int}, H} - \bm{f}_I^{\text{ext}, H} 
 \label{eq:G^H}
 \end{equation}
 and for the analogously defined $G^h$ on the micro domain 
 \begin{equation}
  G^h(\bm{u}^h, \delta \bm{u}^h) \ = \ \sum_{I=1}^{n_{\text{node}}} \left(\delta \bm{d}_I^h \right)^T \bm{f}_I^h
  \quad \text{with} \quad 
  \bm{f}_I^h = \bm{f}_I^{\text{int}, h} - \bm{f}_I^{\text{ext}, h} \, .
  \label{eq:G^h}
 \end{equation}
 The internal and external nodal force vectors are defined for node $I$ of macro element $K$ with volume $|K|$ and surface $\partial K$ and for node $I$ of micro element $T$ with volume $|T|$ and surface $\partial T$ 
 \begin{equation}
 \bm{f}_I^{\text{int}, H} \ = \int_{|K|}  \bm{B}_I^T \bm S^H \, \text{d}V \ , \quad \bm{f}_I^{\text{int}, h} \ = \int_{|T|} \bm{B}_I^T \bm S^h \, \text{d}V
 \label{eq:f_I^int}
 \end{equation}
 
 \begin{equation}
 \bm{f}_I^{\text{ext}, H} \ = \ \int_{\partial K}^{} {N}_I \bm{\overline{t}} \, \text{d}A \ , \quad \bm{f}_I^{\text{ext}, h} \ = \ \int_{\partial T}^{} N_I \bm{t} \, \text{d}A \, .
 \label{eq:f_I^ext}
 \end{equation}
 For the computation of $\bm{f}_I^{\text{int}, H}$ on the macro domain the second Piola-Kirchhoff stress $\bm S^H$ is obtained according to \eqref{eq:macro-mean-stress-S}, hence, along with our notation for the discretized case, by 
 \begin{equation}  
  {\color{black} \bm S^H \ = \ \langle \bm F^{h} \rangle_0^{-1} \langle \bm P^h \rangle_0} \, .
 \label{eq:macro-mean-stress}
 \end{equation}

%-----------------------------------------------------------------------------------------------
\subsection{The micro-to-macro stiffness transfer and assembly}
\label{subsec:micro-2-macro-stiffness-transfer}

Since $G^H$ is nonlinear in $\bm u^H$, it is linearized at the current macro displacement state $\bm{u}^H$  
\begin{equation}
\text{Lin} \left[ G^H(\bm{u}^H, \delta \bm{u}^H) \right] \ = \ {G}^H(\bm{u}^H, \delta \bm{u}^H) +  \text{D}G^H(\bm{u}^H, \delta \bm{u}^H) \cdot \Delta\bm{u}^H \ = \ 0 \, . 
\label{L(GH)}
\end{equation} 
 The first term ${G}^H(\bm{u}^H, \delta \bm{u}^H)$ is the evaluation of $G^H$ for the current macro displacement field, the second term is the linearization of $G^H$ in the direction of the incremental macro displacements $\Delta \bm{u}^H$ ({\color{black} a Directional- or Gateaux-derivative}). \\
 In analogy to the linear case, the key step of nonlinear FE-HMM is to approximate in the numerical quadrature of \eqref{eq:ModifiedBilinearForm-1} the macro values at the quadrature points by the volume-averaged integral of micro scale quantities according to \eqref{eq:ModifiedBilinearForm-3}  
\begin{align}
 \text{D}{G}^H(\bm{u}^H, \delta \bm{u}^H) \cdot \Delta\bm{u}^H  
\ &= 
\sum_{K\in \mathcal T_H} \sum_{l=1}^{N_{qp}} \omega_{K_{\delta_l}} 
\left[ \delta {\bm{E}^H} \colon  {\mathbb{C}_T^{0,h}}  \colon \Delta {\bm{E}^H} + {\bm{S}^H} \colon \Delta \delta {\bm{E}^H} \right] \, \text{d}V 
\label{eq:ModifiedBilinearForm-1} 
\\
\ &\approx 
\sum_{K\in \mathcal T_H} \sum_{l=1}^{N_{qp}} \omega_{K_{\delta_l}} 
\left[ \dfrac{1}{{\color{black}|K_{\delta_l}|}} 
\int_{K_{\delta_l}} \left[\delta \bm{E}^{h} \colon \mathbb{C}^{\epsilon}_T  \colon \Delta \bm{E}^{h} + \bm{S}^{h} \colon \Delta \delta \bm{E}^{h} \right] \right] \text{d}V .
 \label{eq:ModifiedBilinearForm-3}
\end{align}
The approximation reflects that the constitutive law in terms of the free energy function $\psi^{0}(\bm x)$, and the derived quantities of stress $\bm S^H$ and the tangent ${\mathbb{C}_T^{0,h}}$ are not known.\\[5mm]
The micro solution $\bm u^h$ and derived quantities $\bm{E}^{h}$, $\bm{S}^{h}$ are obtained 
\begin{itemize}
	\item on micro sampling domains $K_{\delta_l}=\bm x_{K_{\delta_l}} + \delta \, [-1/2, +1/2]^{n_{\text{dim}}}$, $\delta \geq \epsilon$, volume $|K_{\delta_l}|$
	\item centered at quadrature points $\bm x_{K_{\delta_l}}$ of $K$, $l=1, \ldots,  N_{qp}$ with weights $\omega_{K_{\delta_l}}$,
\end{itemize}
for a visualization see Fig.~\ref{fig:MicMac-Problem-Meshing-etc}. The micro sampling domains render the additive contribution to the tangential stiffness matrix of the macro finite element. 

The approximation in numerical quadrature is consistent with the Hill-Mandel postulate of energy equivalence. It is, however, not sufficient. Since the pointwise macroscopic energy density is replaced by microscopic energy density in a volume, boundary conditions consistent to the postulate must be chosen. 

The macro element tangential stiffness matrix contribution $\bm k^{H}_{T,IJ}$ for nodes $I$ and $J$ of element $K$ are obtained --in analogy to the linear case-- by inserting finite element shape functions into \eqref{eq:ModifiedBilinearForm-1}, hence $\delta \bm u^H(\bm N_I)$ and $\Delta \bm u^h(\bm N_J)$, which results in\footnote{Here and in the following small letters (for stiffness matrices $\bm k$, force vectors $\bm f$, displacements $\bm d$) refer to the element level, capital letters to global, assembled quantities, which is applied for micro scale ($h$) and macro scale quantities ($H$).}
\begin{eqnarray}          
\bm k^{H}_{T,IJ,K} 
&=& \text{D}{G}^{H}(\bm{u}^H, \delta \bm{u}^H (\bm N_J) ) \cdot \Delta\bm{u}^H (\bm N_I)     \\
&=& \label{eq:k-mac-element-5}
\sum_{l=1}^{N_{qp}} \dfrac{\omega_{K_{\delta_l}}}{|K_{\delta_l}|} 
\, \left( \bm D^{h(I)} \right)^T \, \bm K^{h}_{T,K_{\delta_l}} \,  \bm D^{h(J)} \, ,  \\
& & \text{where} \quad \bm D^{h(I)}=\left(\,\bm D^{h(I,x_1)} | \bm D^{h(I,x_2)} | \bm D^{h(I,x_3)} \, \right)  \quad \text{for} \, \, {n_{\text{dim}}}=3 \, ,
\\           
\bm k^{H}_{T,K} &=& \sum_{l=1}^{N_{qp}} \dfrac{\omega_{K_{\delta_l}}}{|K_{\delta_l}|} 
\, \, \bm T^{T}_{K_{\delta_l}}\, \bm K^{h}_{T,K_{\delta_l}} \, \bm T_{K_{\delta_l}} \, , 
\label{eq:k-mac-element-6}  
\\
& & \mbox{where} \quad \bm T_{K_{\delta_l}} = \bigg[ \Big[ \big[ \bm D^{h(I,x_i)} \big]_{i=1,\ldots,n_{\text{dim}}} \Big]_{I=1, \ldots, N_{\text{node}}} \bigg] \label{eq:k-mac-element-7} \,.             
\end{eqnarray} 
The assembly of the global micro stiffness matrix implies the arrangement of $\bm D^{h(I)}$ in different columns for $I=1,\ldots, N_{\text{node}}$, which results in the transformation matrix $\bm T_{K_{\delta_l}}$.

With $M_{\text{mic}}$ for the number of nodes in the micro domain, and $N_{\text{node}}$ for the number of elements on a macro element, it can be easily verified that the matrix dimensions of $\bm T_{K_{\delta_l}} \in \mathbb{R}^{(M_{\text{mic}} \cdot {n_{\text{dim}}}) \times (N_{\text{node}} \cdot {n_{\text{dim}}})}$,  $\bm K^{h}_{T,K_{\delta_l}} \in \mathbb{R}^{(M_{\text{mic}}\cdot {n_{\text{dim}}}) \times (M_{\text{mic}} \cdot {n_{\text{dim}}})}$, and of $\bm k^{H}_{T,K} \in \mathbb{R}^{(N_{\text{node}} \cdot {n_{\text{dim}}})\times (N_{\text{node}} \cdot {n_{\text{dim}}})}$ are consistently used in \eqref{eq:k-mac-element-6}, and that $\bm T_{K_{\delta_l}}$ is an operator of micro-to-macro stiffness transfer and dimensional compression.
 
The single columns of this transformation matrix contain dimensionless displacements of the micro nodes following from macro unit displacement states in terms of shape functions. The superscript $(I,x_i)$ encodes the unit displacement state at macro node $I, I=1,\ldots, N_{\text{node}}$ in the direction of $x_i, i=1, \ldots, n_{\text{dim}}$.
 
The transformation matrix and its derivation provide the information on how to carry out the computations for its entries; first, on the driver of the micro problem in terms of shape functions as described, second, on the corresponding solutions $\bm d^{h (I, x_i)}$, which build up $\bm T_{K_{l}}$. These two characteristics are not sufficient to properly solve the micro problem. The missing links are the coupling conditions, cf.~ Sec.~\ref{subsec:Charactistics-and-conditions-for-CompHomog}.    

 \begin{itemize} 
     \item {\bf Energetically consistent BCs}. They must be applied to the micro domain to comply (the Hill-Mandel postulate of) the micro-macro energy equivalence. The most frequently used BCs are periodic (PBC), linear displacement or kinematically uniform BCs (KUBC/Dirichlet), and constant traction (Neumann) BCs.
     \item  {\bf Order of computational homogenization}. If the macro deformation imparted on the micro domain is homogeneous, the micro problem resembles the unit cell problem of (asymptotic) mathematical homogenization, which is the foundation of FE-HMM. In computational mechanics this setting is frequently referred to as first order computational homogenization. For that aim the macro displacements $\bm u^H$ are linearized at quadrature points $\bm x_{K_{\delta_l}}$ according to  
     \begin{equation}
      \bm u^{H}_{\text{lin}, K_{\delta_l}} = \bm u^H (\bm x_{K_{\delta_l}}) 
       + (\bm x - \bm x_{K_{\delta_l}}) \cdot \nabla \bm u^H(\bm x_{K_{\delta_l}}) \, .
     \label{eq:linearization-uH}
     \end{equation}
 \end{itemize}
With these coupling conditions the embedding of the micro problem into the macro continuum is defined and the micro problem can be solved. The task is cast into a saddle-point problem with the corresponding Lagrange-functional to be minimized.   

For the micro problem driven by the linearized macroscopic unit displacement states  $(I,x_i)$ and thereby subject to energetically consistent BC, solve
\begin{equation}
\mathcal{L} (\bm{D}^{h(I, x_i)},\bm \Lambda^{(I, x_i)}) = \dfrac{1}{2} \bm{D}^{h(I, x_i)T}  \bm{K}_{T}^{h} \bm{D}^{h(I, x_i)}  + \bm \Lambda^{(I, x_i)\,T} \bm G \, \left(\bm D^{h(I, x_i)} - \overline{\bm d}^{H(I, x_i)}\right) \ \rightarrow \ \text{min.}
\label{eq:Lagrange-Functional-for-StiffnessTransfer}
\end{equation}
where\footnote{The basis of the macro element quantity $\overline{\bm d}^H$ is that one of the micro space; therefore and to avoid notational confusion with assembled macro displacements, we employ small letter $\bm d$.} $\overline{\bm d}^{H(I, x_i)}$ is the prescribed nodal displacement vector that follows from linearization in \eqref{eq:linearization-uH}, and where the nodal constraints are encoded in matrix $\bm G$. For convenience we have dropped the index of $\bm K_T^h$ but keep in mind that it refers to a micro problem at macro quadrature point $K_{\delta_l}$. The first variation of $\mathcal{L}$ with respect to both $\bm{D}^{h(I, x_i)}$ and $\bm \Lambda^{(I, x_i)}$ results in the set of linear equations to  be solved for $I=1,...,N_{\text{node}}$, $i=1,..., {n_{\text{dim}}}$: 
\begin{equation}
\left[ \begin{array}{c c}
\bm{K}_{T}^{h} & \bm{G}^T \\
\bm{G} & \bm{0}
\end{array} \right]
\left[ \begin{array}{c}
\bm{D}^{h(I, x_i)} \\
\boldsymbol{\Lambda}^{(I, x_i)}
\end{array} \right] \ = \
\left[ \begin{array}{c}
\bm{0} \\
\bm{G} \, \overline{\bm{d}}^{H(I, x_i)}
\end{array} \right] \, ,
\label{ch_3_eq_30}
\end{equation}
where for PBC condition $\bm{G}^T \bm{\Lambda}^{(I, x_i)}  \ = \ \bm 0$ is fulfilled by the periodic displacements and antiperiodic tractions.  
\\[1mm] 
The solution vector contains the dimensionless micro displacements $\bm{D}^{h(I, x_i)}$, and the Lagrange multipliers $\bm{\Lambda}^{(I, x_i)}$ enforce the coupling condition. With $\bm{D}^{h(I, x_i)}$ the transformation matrix $\bm{T}_{K_{\delta_l}}$ and thereby the macro element stiffness matrix can be built as described by \eqref{eq:k-mac-element-6} and \eqref{eq:k-mac-element-7}. Assembling the macro scale element contributions results in the global system of equations 
\begin{equation}
\Assem \ \sum_{I=1}^{N_{\text{node}}} \sum_{J=1}^{N_{\text{node}}} \left(\delta \bm{d}_I^H \right)^T \left( \bm{k}^{H}_{T,IJ} \Delta \bm{d}^H_J + \bm{f}_I^{H} \right) \ = \ 0 \, .
\label{eq:Assem_k_T_IK_mac} 
\end{equation}
  
For PBC it holds for node $p$ and its counterpart $q$ at opposite edges (for ${n_{\text{dim}}}=2$) and at opposite faces (for ${n_{\text{dim}}}=3$) of a micro domain $K_{\delta_l}$ the constraint condition with ${\bm u}^{H}_{\text{lin}, K_{\delta_l}}$ according to \eqref{eq:linearization-uH} 
\begin{eqnarray} 
\label{eq:ConstraintCondition-phi_h-phi_H}
\left(\bm u^{h}_{K_{\delta_l}} -  {\bm u}^{H}_{\text{lin}, K_{\delta_l}} \right)(p)
&=&
\left(\bm u^{h}_{K_{\delta_l}} -  {\bm u}^{H}_{\text{lin}, K_{\delta_l}} \right)(q) 
\\
\label{eq:ConstraintCondition-GtimesAlphaMinusBeta}
\bm G \, \left(\bm D^{h}_{K_{\delta_l}} - \overline{\bm d}^{H}_{K_{\delta_l}} \right)
&=& \bm 0 \, , 
\end{eqnarray}
where \eqref{eq:ConstraintCondition-GtimesAlphaMinusBeta} is the Lagrange multiplier counterpart, in which matrix $\bm G$ encodes both the macro-micro coupling condition of homogeneous deformation mediated by $\overline{\bm d}^{H}_{K_{\delta_l}}$ and the non-redundant periodic coupling conditions for partner nodes $p$, $q$ on the micro domain boundary. 

Notice that the micro problem and its solution in Sec.~\ref{subsec:micro-2-macro-stiffness-transfer} merely serves the purpose of building the stiffness carrier $\bm T_K$. For the solution of the macro problem (i) macro stress $\bm S^H$ from averaged micro stresses and (ii) macro stiffness $\bm k^H_T$ from $\bm K^h_{T}$ must be available, namely for the micro state fulfilling the weak form of balance of linear momentum. For that purpose the micro problem is to be solved which is described in the following  Sec.~\ref{subsec:microprob-solution}.  
%----------------------------------------------------------------------------------------------------- 
\subsection{The microproblem and its solution}    
\label{subsec:microprob-solution} 
The micro problem is solved on the RVE of each macro quadrature point $K_{\delta_l}$, $l=1, ..., n_{qp}$ in all finite elements $K$ of the macro triangulation ${\cal T}^H$.  

The solution of the nonlinear system of equations on the micro domain from \eqref{eq:G^h} requires the linearization at the current micro displacement state $\bm{u}^h$ 
 \begin{equation}
\text{Lin} \left[ G^h(\bm{u}^h, \delta \bm{u}^h) \right] \ = \ {G}^h(\bm{u}^h, \delta \bm{u}^h) + \text{D}{G}^h(\bm{u}^h, \delta \bm{u}^h) \cdot \Delta\bm{u}^h \ = \ 0 \, , 
 \label{ch_3_eq_08}
 \end{equation}
 with the definitions already introduced in the context of the linearization of $G^H$.
 
The approximation of $\text{D}{G} \cdot \Delta \bm{u}$ is given by  
 \begin{equation}
 \text{D}{G}^h(\bm{u}^h, \delta \bm{u}^h) \cdot \Delta\bm{u}^h \ = \ \int_{T} \left( \delta \bm E^h \colon \mathbb{C}_T^\varepsilon \colon \Delta \bm{E}^h + \bm{S}^h \colon \Delta\delta\bm{E}^h \right) \, \text{d}V \, .
 \label{ch_3_eq_09}
 \end{equation}
 Inserting the approximations into the framework of the finite element method leads to   
 \begin{equation}
  \text{D}{G}^h(\bm{u}^h, \delta \bm{u}^h) \cdot \Delta\bm{u}^h \ = \ \sum_{I=1}^{n_{\text{node}}} \sum_{K=1}^{n_{\text{node}}} \left(\delta \bm{d}_I^h \right)^T \underbrace{\int_{T}^{} \left( \bm{B}_I^T \mathbb{C}_T^\varepsilon \bm{B}_K + \hat{\bm{G}}_{IK} \right) \, \text{d}V}_{\displaystyle =:\bm{k}_{T,IK}^{h}} \Delta\bm{d}^h_K \, , 
\label{eq:DG^h}
 \end{equation}
 with the tangential element stiffness matrix contribution of nodes $I$ and $K$ and the initial stress matrix $\hat{\bm{G}}_{IK}$ reflecting the contribution of the current stress state to the stiffness of a deformed structure  
 \begin{equation}
 \hat{\bm{G}}_{IK} \ = \ \hat{S}_{IK}^h \bm{1} \qquad \text{with}
 \label{eq:G_IK-InitStressMatr}
 \end{equation}
 \begin{align}
 \nonumber
 \hat{S}^h_{IK} \ = \ & S_{11}^h N_{I,1}N_{K,1} + S_{22}^h N_{I,2}N_{K,2} + S_{33}^h N_{I,3}N_{K,3} + S_{12}^h (N_{I,1}N_{K,2} + N_{I,2}N_{K,1}) + \\
 & S_{13}^h (N_{I,1}N_{K,3} + N_{I,3}N_{K,1}) + S_{23}^h (N_{I,2}N_{K,3} + N_{I,3}N_{K,2}) \, .
 \label{eq:S_IK}
 \end{align}
 The corresponding global finite element approximation on the micro scale is obtained by assembling the  contributions of all micro finite elements (number: $\text{num}_{ele}$)
 \begin{equation}
 \Assem \ \sum_{I=1}^{n_{\text{node}}} \sum_{K=1}^{n_{\text{node}}} \left(\delta \bm{d}_I^h \right)^T \left( \bm{k}^{h}_{T,IK} \Delta \bm{d}^h_K + \bm{f}_I^{h} \right) \ = \ 0 \, .
 \label{eq:Assem_k_T_IK_mic} 
 \end{equation}  
 
The assembled system of equations \eqref{eq:Assem_k_T_IK_mic} is subject to the coupling conditions, hence 
\begin{equation}
%\left(\delta \bm{D}^h\right)^T 
%\left( 
\bm{K}_{T}^{h} \, \Delta \bm{D}^h + \bm{F}^{h} 
%\right) 
\ = \ \bm 0 \qquad \text{subject to} \qquad \bm{G} \left(\Delta \bm{D}^h - \overline{\bm d}^H\right)=\bm 0 \, .
\label{eq:set-linear-eq-microscale}
\end{equation}
In the incremental setting of linearized equations along with the constraint matrix $\bm G$, the constraint for energetically consistent BCs is written in terms of the yet unknown micro displacement increment $\Delta \bm{D}^h$ and the macroscopic displacements $\overline{\bm d}^H$ imposed on the RVE boundary nodes to realize the required homogeneous deformation.
   
For the solution of \eqref{eq:set-linear-eq-microscale} a Lagrange-functional is constructed which is to be minimized.
\begin{equation}
 \mathcal{L} (\Delta \bm{D}^h,\Delta \bm \Lambda^h) = \dfrac{1}{2} \Delta \bm{D}^{h\,T} \bm{K}_{T}^{h} \Delta \bm{D}^h 
 + \bm{F}^{h\,T} \, \Delta \bm{D}^h 
 + \Delta \bm \Lambda^{h\,T} \, \bm G \, \left(\Delta \bm D^h - \overline{\bm d}^H\right)   \  \longrightarrow \  \text{min.} 
 \label{eq:Lagrange-functional-for-Equilibrium-Iteration}
\end{equation}

The partial derivative of $\mathcal{L}$ with respect to $\Delta \bm{D}^h$ and with respect to $\Delta \bm{\Lambda}^h$ yield the stationarity conditions
 \begin{equation}
 \left[ \begin{array}{c c}
 \bm{K}_{T}^{h} & \bm{G}^T \\
 \bm{G} & \bm{0}
 \end{array} \right]
 \left[ \begin{array}{c}
 \Delta \bm{D}^{h} \\
 \Delta \bm{\Lambda}^{h}
 \end{array} \right] \ = \
 \left[ \begin{array}{c}
 -\bm{F}^{h} \\
 \bm{G} \, \overline{\bm{d}}^{H}
 \end{array} \right] \, .
 \label{eq:Saddlepoint-problem-Set-of-eqs}
 \end{equation}
 Therein, the first set of equations
 \begin{equation}
 \bm{K}_{T}^{h} \Delta \bm{D}^h = 
 \underbrace{- \bm{F}^{h} - \bm{G}^T \Delta \bm{\Lambda}^h}_{\displaystyle =: - \bm R^h}
  \label{ch_3_eq_19-1}
 \end{equation}
 corresponds to \eqref{eq:set-linear-eq-microscale}, the term $\bm{G}^T \Delta \bm{\Lambda}^h$ describes the change of the external node forces (Lagrange multipliers enforcing the coupling), which, together with the internal force vector $\bm{F}^{h}$, render the residual $\bm R^h$.
 In the case of PBC the change of these external node forces is not considered in the linearization, since $\bm{G}^T \Delta \bm{\Lambda}^h=\bm 0$ is fulfilled by periodic displacements and antiperiodic tractions. 
  
The solution of \eqref{eq:Saddlepoint-problem-Set-of-eqs} yields the updates
 \begin{align}
 \label{ch_3_eq_22}
 \bm{D}^h \ &\longleftarrow \ \bm{D}^h + \Delta \bm{D}^h \, , \\
 \bm{\Lambda}^h \ &\longleftarrow \ \bm{\Lambda}^h + \Delta \bm{\Lambda}^h \,.
 \label{ch_3_eq_23}
 \end{align}
 Iterations along with updates stop for sufficiently small residuals. Then the computation of the micro-to-macro stiffness transfer by means of the current tangential stiffness matrix and by means of the current stress transfer is carried out.  
   
The application of energetically consistent BCs fulfilling the Hill-Mandel condition in the FE-HMM framework are analyzed in \cite{EidelFischer2018a} for linear elasticity along with convergence and error analysis, for finite deformations see \cite{Saeb-Steinmann-Javili-2016}. 

\bigskip 

{\bf Remarks} 
\\[1mm]
The differences in the Lagrange functionals of \eqref{eq:Lagrange-Functional-for-StiffnessTransfer} in comparison to \eqref{eq:Lagrange-functional-for-Equilibrium-Iteration} reflects their different purpose; \eqref{eq:Lagrange-functional-for-Equilibrium-Iteration} is constructed for the continuation of a Newton iteration by linearized equations in \eqref{eq:Saddlepoint-problem-Set-of-eqs}, whereas \eqref{eq:Lagrange-Functional-for-StiffnessTransfer} is constructed from truely linear equations; for a state fulfilling the balance equations after having solved the nonlinear set of algebraic equations, it is now purely the task of the micro-to-macro stiffness transfer. The stiffness transfer coincides with that one in linear homogenization. For the linear case the tangential stiffness matrix boils down to the deformation- and stress-independent  stiffness matrix, \cite{EidelFischer2016}, \cite{EidelFischer2018}. 
\\[1mm] 
The FE$^2$ method realizes the micro-macro stiffness transfer in terms of the homogenized elasticity tensor/tangent moduli, which are obtained by means of the micro stiffness matrix. FE-HMM in contrast directly draws on the micro stiffness matrices as described. 
 
%-----------------------------------------------------------------------------------------------
\section{A priori error estimates}
%---------------------------------------------------------------------------------------------
  \label{sec:Apriori-Estimates}
 
By virtue of its roots in mathematical homogenization, the FE-HMM exhibits a priori estimates; for the elliptic case in \cite{E-Ming-Zhang-2005}, \cite{Ohlberger2005}, for the elliptic case of linear elasticity in a geometrical linear setting to \cite{Assyr2006}, \cite{Assyr2009}. 
  
The total FE-HMM error can be decomposed into three parts
  \begin{equation}
  || \bm u^0 - \bm u^H || \, \leq \, \underbrace{|| \bm u^0 - \bm u^{0,H} ||}_{\displaystyle e_{\text{mac}}} 
  \, + \, \underbrace{|| \bm u^{0,H} - \widetilde{\bm u}^H ||}_{\displaystyle e_{\text{\text{mod}}}}
  \, + \, \underbrace{|| \widetilde{\bm u}^H - \bm u^H ||}_{\displaystyle e_{\text{\text{mic}}}} \, ,
  \label{eq:Error-decomposition-mac-mod-mic}
  \end{equation}
  where $e_{\text{mac}}$, $e_{\text{mod}}$, $e_{\text{mic}}$ are the macro error, the modeling error, and the micro error.
  
  Here, $\bm u^0$ is the solution of the homogenized problem, $\bm u^H$ is the FE-HMM solution, $\bm u^{0,H}$ is the standard (single-scale) FEM solution that is obtained through exact $\mathbb{C}^{0}$ (in linear elasticity); and $\widetilde{\bm u}^H$ is the FE-HMM solution obtained through exact micro solutions.
  
  For sufficiently regular problems the following a priori estimates hold in the $L^2$-norm, the $H^1$-norm and the energy-norm:  
  \begin{eqnarray}
  || \bm u^0 - \bm u^H ||_{L^2(\mathcal{B}_0)} &\leq& C\left( H^{p+1} + \left(\dfrac{h}{\epsilon}\right)^{2q} \right) + e_{\text{mod}}    \, ,
  \label{eq:Total-Error-estimate-L2} \\
  || \bm u^0 - \bm u^H ||_{H^1(\mathcal{B}_0)} &\leq& C\left( H^p + \left(\dfrac{h}{\epsilon}\right)^{2q} \right) + e_{\text{mod}}   \, ,
  \label{eq:Total-Error-estimate-H1} \\
  || \bm u^0 - \bm u^H ||_{A(\mathcal{B}_0)} &\leq& C\left( H^p + \left(\dfrac{h}{\epsilon}\right)^{2q} \right) + e_{\text{mod}}   \, .
  \label{eq:Total-Error-estimate-Energy}  
  \end{eqnarray}
  
  For $e_{\text{mod}}$ in \eqref{eq:Total-Error-estimate-L2}--\eqref{eq:Total-Error-estimate-Energy} it holds 
  
  \begin{equation}
  \label{eq:ModelingError}
  e_{\text{mod}} = \left\{ \begin{array}{ll}
  0 & \mbox{for periodic coupling with} \,\, \delta/\epsilon \in \mathbb{N} \\
  {\color{black}C\,\dfrac{\epsilon}{\delta}}  & \mbox{for Dirichlet coupling with} \,\, \delta > \epsilon \end{array} \right. \,    
  \end{equation} 
  {\color{black} given that the hypotheses hold, that the elasticity tensor $\mathbb{C}^{\epsilon}$ is periodic on the RVE and, that the micro solution is sufficiently smooth, \cite{JeckerAbdulle2016}.}
    
  The modeling error for Dirichlet coupling in \eqref{eq:Total-Error-estimate-L2}--\eqref{eq:Total-Error-estimate-Energy} is due to boundary layers \cite{E-Ming-Zhang-2005} (Thm. 1.2), \cite{Assyr2009}. So even for $H \rightarrow 0$ and $h \rightarrow 0$ there is a residual error. 

The second terms in the estimates  \eqref{eq:Total-Error-estimate-L2}--\eqref{eq:Total-Error-estimate-Energy} each represent the micro error which is propagated to the macro scale. The deviations from the convergence order of the micro error on the micro scale is threefold; first, in that the error in the $L^2$-norm is of order $2q$ instead of $q+1$, second, in that errors in the $H^1$- and energy-norm are of the order of $2q$ instead of $q$ and therefore and third, in that it coincides for the $L^2$-norm with the $H^1$-/energy-norm.

For the particular case of linear elasticity it holds
\begin{equation}
    \label{eq:Estimate-ElasticityTensor}
    || \mathbb{C}_{ijkl}^0 - \mathbb{C}_{ijkl}^{0,h} || \leq C  \left(\dfrac{h}{\epsilon}\right)^{2q}    
\end{equation}
with the exact homogenized elasticity tensor $\mathbb{C}_{ijkl}^0$ and its approximation $\mathbb{C}_{ijkl}^{0,h}$ obtained on a grid with element size $h$.  

The unified estimates covering both the macro error as well as the micro error enable uniform micro-macro ($h$-/$H$-) mesh refinements for optimal convergence order but for minimal computational costs. 
\begin{Table}[htbp]
	\begin{minipage}{16.5cm}  
		\centering
		\renewcommand{\arraystretch}{1.2} 
		\begin{tabular}{c|cc}
			\hline 
			macro, micro FEM  &  $L^2$-norm                    & $H^1$-/energy-norm \\
			\hline  
			$P^p$, $P^q$      & $N_{\text{mic}} = (N_{\text{mac}})^{(p+1)/2q}$ & $N_{\text{\text{mic}}} = (N_{\text{mac}})^{p/2q}$    \\
			\hline                 
		\end{tabular}
		\newline 
	\end{minipage}
	\caption{{\color{black}{\bf} Optimal uniform micro-macro refinement strategies: full order for minimal effort. {\color{black} $N_{\text{mic}}$ denotes the number of unknowns on the micro scale, $N_{\text{mac}}$ on the macro scale.}
			\label{tab:Best-mic-mac-RefinementStrategies}}}
\end{Table}

Table \ref{tab:Best-mic-mac-RefinementStrategies} displays the optimal uniform micro-macro refinement strategies for the error in the  $L^2$-norm and the $H^1$-/energy-norm. Of course, the strategies' dependency on the polynomial orders of the macro and micro shape functions $p$ and $q$ rely on full regularity of the corresponding BVPs. In the present context of nonlinear homogenization, the additional question arises, whether the estimates apply at all. For that aim we measure the convergence and compare it with the fully linear case in order to properly identify the source of possible order reductions. 

For linear problems in elasticity the estimates  \eqref{eq:Total-Error-estimate-L2}--\eqref{eq:ModelingError}, \eqref{eq:Estimate-ElasticityTensor} as well as the optimal refinement strategies according to Tab.~\ref{tab:Best-mic-mac-RefinementStrategies} were assessed and for sufficient regularity confirmed in \cite{EidelFischer2016}, \cite{JeckerAbdulle2016}, \cite{EidelFischer2018}, \cite{EidelFischer2018a}. 
 
\input{AlgorithmicSolutionConcepts}

\input{Examples}

\section{Summary} 

The main results of the present paper shall be summarized.
 
\begin{enumerate}
 \item A Finite Element Heterogeneous Multiscale Method FE-HMM was developed for geometrical nonlinearity and hyperelastic solids in a Lagrangean setting.    
 \item Existing a priori estimates for the fully linear case (linear elasticity in a geometrical linear frame) were assessed for the present nonlinear case in numerical tests and confirmed for sufficient regularity of the BVPs. In either case the convergence results in the nonlinear elastic regime are close to those of the fully linear case for the investigated examples.
 \item The staggered, micro-macro solution concept in terms of a two-level Newton algorithm was revised along with the proposal of an alternative scheme.
       \begin{itemize}
        \item The standard approach in terms of nested loops can safely and efficiently be replaced by direct alternations between micro and macro solution iterations. The fully converged solution of micro problems for each macro iteration is neither necessary for achieving convergence nor favorable for the efficiency.
        \item For the considered problems the novel concept (i) is robust with respect to  large loading steps, (ii) exhibits speedup factors in the range of 1.05 up to 2.6 in 2D, and in the range from 1.7 to 2.0 in 3D, where the speedup is on average higher in 3D than in 2D. 
       \item The observed minimal speedup is above 1.0, such that the novel concept does not slow down simulations. Since the speedup of 1.0 is clearly not a granted lower bound, monitoring the macro convergence on-the-fly in terms of the required iterations can avoid deficits in efficiency. 
       \item A particular appeal of the proposed concept is its simplicity and ease of implementation. 
       \end{itemize}  
 \item The advantages of the novel concept was demonstrated for FE-HMM simulations. It equally applies for FE$^2$ for the equality of the methods despite a minor difference in micro-macro stiffness transfer.
\end{enumerate} 
It is expected that the novel speedup concept can be used in different applications as e.g. in coupled problems such as electro-/magneto-/mechanical-coupling, or generally, in any other nonlinear problem set for which a staggered, two-level algorithmic solution concept in terms of an embedded loop is standard. 

Inelastic problems in homogenization require a separate analysis in order to arrive at valid conclusions on the applicability of the proposed speedup concept, which will be presented elsewhere.  

\bigskip

{\bf Acknowledgments.} Bernhard Eidel acknowledges support by the Deutsche Forschungsgemeinschaft (DFG) within the Heisenberg program (grant no. EI 453/2-1). Simulations were performed with computing resources granted by RWTH Aachen University under project ID prep0005. This has enabled the present work.

\bibliography{literature_homogenization}

%% file: Frontpage-1.tex
\vspace{-3mm}
\begin{center}
  {\bf \large A Nonlinear Finite Element Heterogeneous Multiscale Method} \\[2mm]
  {\bf \large for the Homogenization of Hyperelastic Solids and}\\[2mm]
  {\bf \large a Novel Staggered Two-Scale Solution Algorithm} 
\end{center}

\vspace{4mm}
\ce{Bernhard Eidel, Andreas Fischer, Ajinkya Gote
}
 
\vspace{4mm}
 
\ce{\small Heisenberg Group, Institute of Mechanics, Department Mechanical Engineering}
\ce{\small Universit\"at Siegen, 57068 Siegen, Paul-Bonatz-Str. 9-11, Germany}
\ce{\small e-mail: bernhard.eidel@uni-siegen.de, phone: +49 271 740 2224, fax: +49 271 740 2436}
\vspace{2mm}

\bigskip

\begin{center}
{\bf \large Abstract}

\bigskip

{\footnotesize
\begin{minipage}{14.5cm}
\noindent
In this paper we address three aspects of nonlinear computational homogenization of elastic solids by two-scale finite element methods. First, we present a nonlinear formulation of the finite element heterogeneous multiscale method FE-HMM in a Lagrangean formulation that covers geometrical nonlinearity and, more general, hyperelasticity. Second, a-priori estimates of FE-HMM, which exist so far only for the fully linear elastic case in solid mechanics, are assessed in the regime of nonlinear elasticity. 
The measured convergence rates agree fairly well with those of the fully linear regime. Third, we revise the standard solution algorithm of FE$^2$ which is a staggered scheme in terms of a nested loop embedding the full solution of the micro problem into one macro solution iteration step. We demonstrate that suchlike staggered scheme, which is typically realized by a nested two-level Newton algorithm, can safely and efficiently be replaced by direct alternations between micro and macro iterations. The novel algorithmic structure is exemplarily detailed for the proposed nonlinear FE-HMM, its efficiency is substantiated by a considerable speedup in numerical tests. 
\end{minipage}
}
\end{center}
 
{\bf Keywords:}
Heterogeneous multiscale method; Finite element method; Macro-to-micro modeling; Homogenization
%\hspace*{6.3cm} 
%\\[2mm]
%vers.\,\today \, at \currenttime\\

%% file: AlgorithmicSolutionConcepts.tex
\section{Novel Algorithmic Solution Concept}
\label{sec:NovelAlgorithm}
 
\subsection{The standard: nested loops}
\label{subsec:standard}
  
\begin{Figure}[htbp]
   \begin{minipage}{16.5cm} 
   \centering 
   \includegraphics[width=15.0cm, angle=0, clip=]{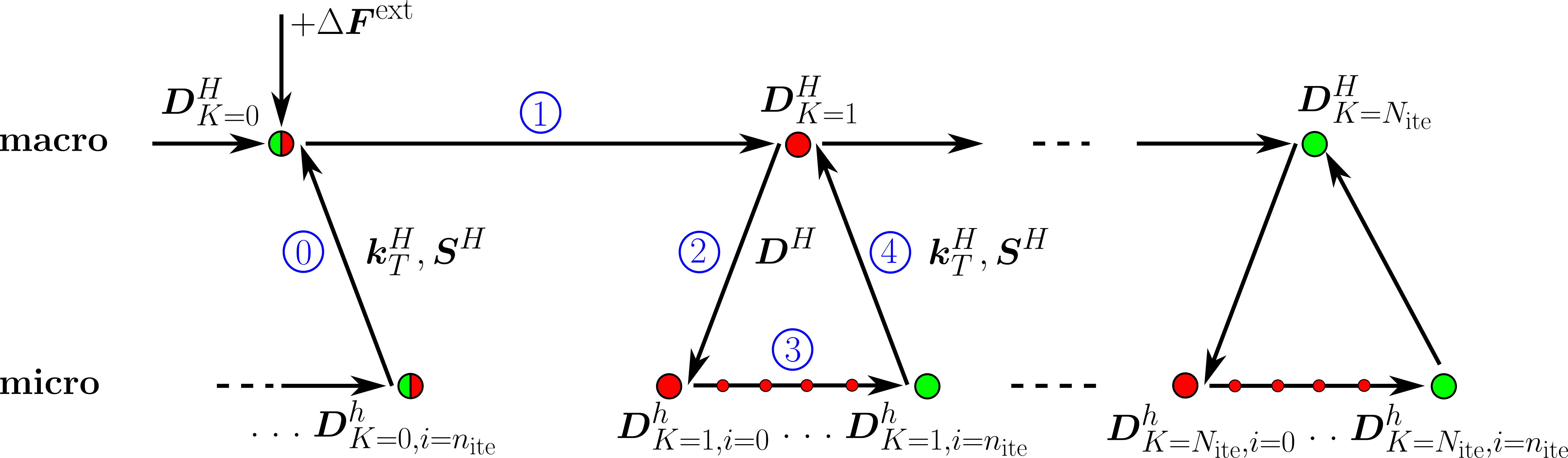}
      \vspace*{2mm} 
   \end{minipage}
 \caption{The conventional staggered scheme based on nested, macro-micro solution loops for the nonlinear homogenization by a two-level FEM. Green (red) circles indicate stages (not) in equilibrium. 
\label{fig:MicMac-Staggered-Solution-scheme-nested}}   
\end{Figure}

The algorithmic standard of nonlinear FE$^2$ is a staggered solution scheme in terms of a nested loop\footnote{In \cite{FeyelChaboche2000}, p.313, referred to as ''interleaved finite element algorithms''} embedding the complete solution of the micro problem into one macro solution iteration step as displayed in pseudocode ''Algorithm 1'' in Sec.~\ref{sec:intro}. The scheme is typically realized by a nested two-level Newton method.

This approach is depicted for a Lagrangean formulation in Fig.~\ref{fig:MicMac-Staggered-Solution-scheme-nested} and goes like this; a loading increment --here in terms of an external force increment $\Delta \bm F^{\text{ext}}$-- shifts both the macro scale as well as the micro scale out of balance of linear momentum, the lights go from green to red.
For the solution of the macro problem, the current states of macro stiffness $\bm k_T^{H}$ and stress $\bm S^H$ are computed on the micro scale and transferred from micro to macro; it is an initialization for the start of a simulation, or an update for the continuation of a simulation, in either case denoted by $(0)$ in  Fig.~\ref{fig:MicMac-Staggered-Solution-scheme-nested}.
The updated macro stiffness and residual --stress enters the internal force vector being part of the residual-- enable the first iteration of the macro Newton, indicated by step $(1)$ in Fig.~\ref{fig:MicMac-Staggered-Solution-scheme-nested} and its resultant macro displacement vector $\bm D^{H}_{K=1}$.
  
In parallel, the fluctuations (for periodic BCs on the RVE) of the converged micro solution $\bm D^{h}_{K=0,i=n_{\text{ite}}}$ from the finished previous loading step are maintained. They are superimposed to the updated displacement field --step (2) in Fig.~\ref{fig:MicMac-Staggered-Solution-scheme-nested}-- to provide favorable start values for solving the micro problems by Newton's method in the current loading step. 
The iterations in the solution of the micro problem are indicated by the counter $i=0,\ldots,n_{\text{ite}}$, step (3) in 
Fig.~\ref{fig:MicMac-Staggered-Solution-scheme-nested}. Despite its finally vanished micro residuals, the obtained micro solutions each represent a pseudo-equilibrium state\footnote{They are notwithstanding in green color in Fig.~\ref{fig:MicMac-Staggered-Solution-scheme-nested} even if the macro Newton is not yet converged.} for the macro residual still being non-zero.

Based on the obtained micro solution, macro (element) stiffness $\bm k_T^{H}$ and stress $\bm S^H$ are passed over to the macro level, step (4) in Fig.~\ref{fig:MicMac-Staggered-Solution-scheme-nested}, which completes the cycle of one macro Newton iteration. These cycles are repeated until the macro Newton is converged, which is indicated by $K=n_{\text{ite}}$. An accuracy polish by a final iteration on the micro scale completes the loading step on both scales. 

The inherent hierarchy in the standard staggered scheme in that the outer loop for the macro Newton can only proceed with an update of micro scale quantities follows from the very multiscale setup with its scale-dependencies; the continuation of the macro solution process after one iteration requires the update of stiffness and stress from the micro scale. It is, however, not clear, whether the macro Newton requires fully converged micro solutions for continuation. 

\subsection{The novel: direct alternations}
\label{subsec:novel}

\begin{Figure}[htbp]
%	\hspace*{5mm}
	\begin{minipage}{16.5cm}
		\centering  
\includegraphics[width=15.0cm, angle=0, clip=]{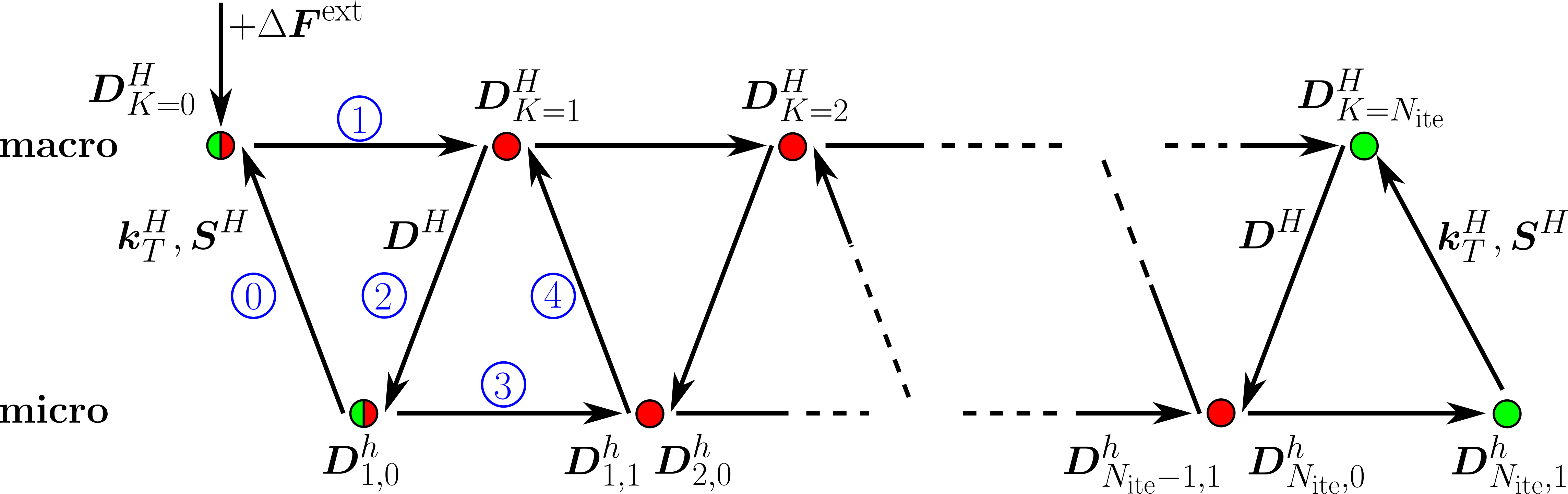}
 		\vspace*{2mm}   
	\end{minipage}
	\caption{The novel staggered scheme with direct macro-micro alternations in the two-level Newton iterations for the nonlinear homogenization 
		by a two-level FEM. Green (red) circles indicate stages (not) fully converged. 
		\label{fig:MicMac-Staggered-Solution-scheme-novel}}   
\end{Figure}

The present work analyzes, whether the nested solution algorithm as described is necessary or at least favorable for the homogenization of hyperelastic solids. The criterion is of purely numerical nature, related to convergence in a qualitative and a quantitative manner. The only physical requirement is that balance laws must be fulfilled both on the micro scale and the macro scale at the end of each load step,  here the balance of linear momentum.

We hypothesize that the continuation of the macro Newton does neither require nor take profit in its overall performance from micro solutions that are fully converged. In that case the standard is an expensive luxury. 

In line with this hypothesis we advocate an alternative route for a novel staggered scheme as depicted in pseudocode ''Algorithm 2'' in Sec.~\ref{sec:intro} and --applied to the mechanical problem-- in  Fig.~\ref{fig:MicMac-Staggered-Solution-scheme-novel}, which promises considerable computational savings.    

The idea is to replace the embedded-loop solution scheme by a loop of direct micro-macro alternations; after each micro Newton iteration the current stiffness $\bm k_T^{H}$ and stress $\bm S^H$ are passed over to the macro scale which is followed --as in the standard-- by a macro iteration along with a displacement update for the micro scale.
  
In the novel version the number of consecutive iterations on the micro level is reduced to one, e.g. the step from $\bm D^h_{1,0}$ to $\bm D^h_{1,1}$ (first index $K$, second index $i$) denoted by (3) in Fig.~\ref{fig:MicMac-Staggered-Solution-scheme-novel}, the rest $i=2,\ldots, n_{\text{ite}}$ of the standard scheme is discarded. After the consecutive macro iteration an updated displacement field yields micro displacements $\bm D^h_{2,0}$ which contain for PBC the fluctuations stored from the state $\bm D^h_{1,1}$. Consequently, the last micro iteration carries out a step from 
$\bm D^h_{N_{\text{ite}},0}$ to $\bm D^h_{N_{\text{ite}},1}$. For an already converged macro solution, an additional accuracy polish of macro stress is realized by this final micro iteration, which is equally carried out in the standard solution algorithm.

%% file: Examples.tex
%---------------------------------------------------------------------------------------------------------
\section{Numerical Examples}
\label{sec:NumericalExamples}
%---------------------------------------------------------------------------------------------------------
  
In this section the nonlinear FE-HMM formulation is assessed with a focus on two aspects.
\begin{enumerate}
 \item[(A)] For the two-level Newton algorithm we compare the performance of the standard version with the novel version using direct micro-macro alternations. 
 \item[(B)] Assessment of the a priori error estimates \eqref{eq:Total-Error-estimate-L2}--\eqref{eq:Total-Error-estimate-Energy} for the case of geometrical nonlinearity along with either linear elasticity or (nonlinear) hyperelasticity. Doing so we keep in mind that the estimates are proved for linear elasticity in a geometrical linear frame. A comparison with the convergence of the fully linear case is necessary to identify, whether a deviation from the theoretical order is due to the nonlinearity of the considered problem or due to the reduced regularity of the micro problem setting itself. Errors are calculated by means of an accurate reference solution obtained on very fine meshes.  
\end{enumerate}  

\begin{Figure}[htbp]
	\begin{minipage}{16.5cm}  
		\centering                                                 
		\includegraphics[width=10cm, angle=0, clip=]{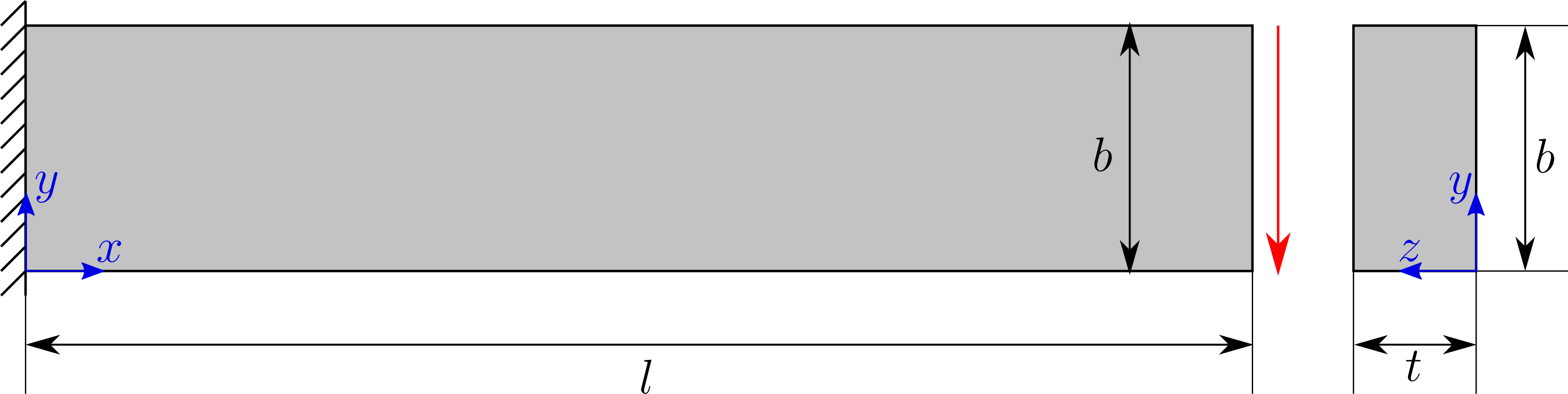}
	\end{minipage}
	\caption{{\bf Macro problem}: Cantilever beam, geometry and boundary conditions. For notational convenience we use the Cartesian coordinates $(x,y,z)$ replacing 
		the previous $(x_1,x_2,x_3)$ coordinate system. 
		\label{fig:macroproblem}}
\end{Figure} 
The macro problem common to all micro problems is a bending-dominated cantilever beam. It exhibits length $l$ in $x$-direction, height $b$ in $y$-direction, and thickness $t$ in $z$-direction. It is clamped at $x$=$0$ and subject to force control by an external line load $\bm f$ or to displacement control at $x$=$l$. If not otherwise stated, it holds $l=5000$~mm, $b=1000$~mm, $t=100$~mm. The RVE side length $\epsilon$ as well as the type and magnitude of loading is defined in each example separately. 

For hyperelasticity we consider an isotropic neo-Hookean constitutive law frequently used for rubber-like materials. The free energy function exhibits the form

\begin{equation}
\label{eq:NH-free-energy}
\psi = \frac{\lambda }{4}\left ( J^{2} - 1 \right )-\left ( \frac{\lambda }{2} + \mu  \right )  \text{ln} \, J + \frac{\mu }{2}\left ( \text{tr} \, \bm C -3 \right ) \, ,
\end{equation}

where $\lambda$ and $\mu$ are the Lam\'{e} parameters, $J=\text{det} \, \bm F$ is the Jacobian and $\bm C$ is the right Cauchy-Green deformation tensor. From $\psi$ and the hyperelasticity relations for stress and the tangent, $\bm S = 2\partial_{\bm C} \psi$, $\mathbb{C}=4 \partial^{2}_{\bm C \, \bm C} \psi$, we obtain
\begin{eqnarray}
\boldsymbol{S} &=& \frac{\lambda }{2}\left ( J^{2} -1\right )\boldsymbol{C}^{-1} + \mu \left ( \boldsymbol{I} - \boldsymbol{C}^{-1} \right ) \, , \\
\mathbb{C}_{ijkl} &=& \left [ \lambda\left ( J^{2}-1 \right ) - 2\mu  \right ]\mathbb{A}_{ijkl} + \lambda J^{2}C^{-1}_{ij}C^{-1}_{kl} \, ,\\
\text{with} \quad 
\mathbb{A}_{ijkl} &=& \frac{\partial \boldsymbol{C}^{-1}}{\partial \boldsymbol{C}} = -\frac{1}{2}\left ( C^{-1}_{ik}C^{-1}_{jl} + C^{-1}_{jk}C^{-1}_{il} \right ) \, . 
\label{eq:Aijkl-NeoHook}
\end{eqnarray}

Notice that the constitutive law is defined on the micro scale; for notational convenience we discard in \eqref{eq:NH-free-energy}--\eqref{eq:Aijkl-NeoHook} the explicit micro scale label of $\epsilon$.  

The considered microstructures consist of two different elastic phases. In view of their considerable differences in morphology, the elasticities of the phases are comprehensively chosen to be same for all examples; the material parameters are given in Tab.~\ref{tab:Material-param-pegasus} for linear elasticity and neo-Hookean hyperelasticity.

\begin{Table}[htbp] 
	\centering
	\begin{tabular}{cccccc}
		\hline
		& \multicolumn{2}{c}{Linear Elasticity}& \ &   \multicolumn{2}{c}{Neo-Hookean}  \\ \hline
		&      E      & $\nu$  & \ &  $\lambda$  &   $\mu$    \\
		phase	& (N/mm$^2$)  &        & \ &  (N/mm$^2$) & (N/mm$^2$) \\ \hline
		1	    &   100\,000    &  0.2   & \ &   27\,777.78  & 41\,666.67   \\
		2       &    40\,000    &  0.2   & \ &   11\,111.11  & 16\,666.67   \\ \hline
	\end{tabular}
	\caption{Material parameters for two-phase microstructures.} 
	\label{tab:Material-param-pegasus}
\end{Table}

The simulations are carried out by a C++ code using MPI-parallelization.

%----------------------------- 
\subsection{Escher's Pegasus} 
\label{subsec:Escher_pegasus}
%----------------------------- 

The first example considers on the macro scale a cantilever beam as displayed in Fig.~\ref{fig:macroproblem}, which is subject at $x=l$ to a lineload of magnitude $200$ N/mm in negative $y$-direction. The material is a periodic, two-phase microstructure as displayed in Fig.~\ref{fig:Escher-Pegasus} similar to the tessellation ''Pegasus N$^{\circ}$105'' of the Dutch graphic artist M.C. Escher. The material's elastic parameters are given in Tab.~\ref{tab:Material-param-pegasus}.

\begin{Figure}[htbp]
   \begin{minipage}{16.0cm}  
   \centering  
         \includegraphics[height=5cm, angle=0, clip=]{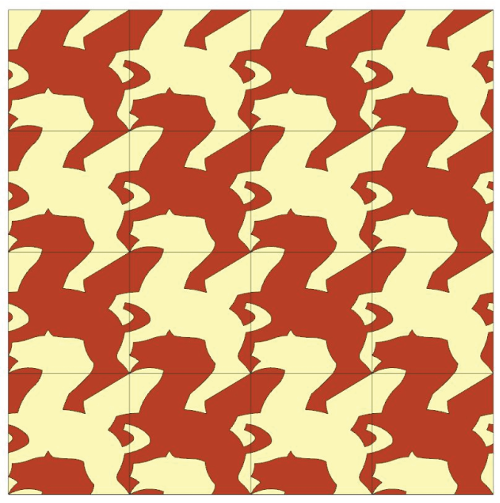}  
         \hspace*{4mm}  
         \includegraphics[height=4.9cm, angle=0, clip=]{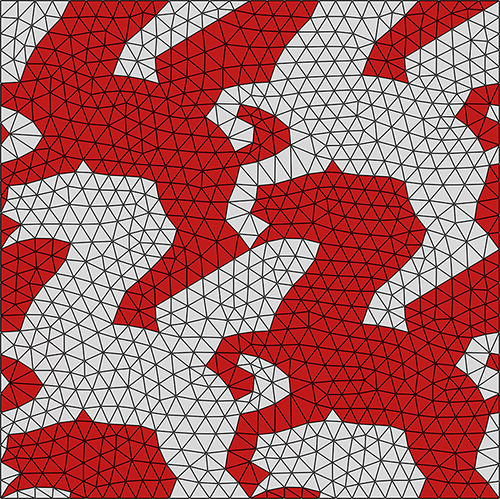}  
   \end{minipage}
\caption{{\bf Escher's Pegasus}: (left) Periodic tessellation with 2$\times$2 unit cells and (right) a triangulation of a unit cell.}  
\label{fig:Escher-Pegasus} 
\end{Figure}
   
\subsubsection{Micro convergence for linear elasticity and Neo-Hookean hyperelasticity} 
%------------------------------------------ 
The reference solution for micro convergence analysis is calculated with ndof=24 on the macro scale domain (5$\times$1 macro elements) and ndof=2\,262\,584 on the micro scale with $\epsilon=110$~mm.  

\begin{Figure}[htbp]
	\centering
	\includegraphics[height=4.5cm, angle=0, clip=]{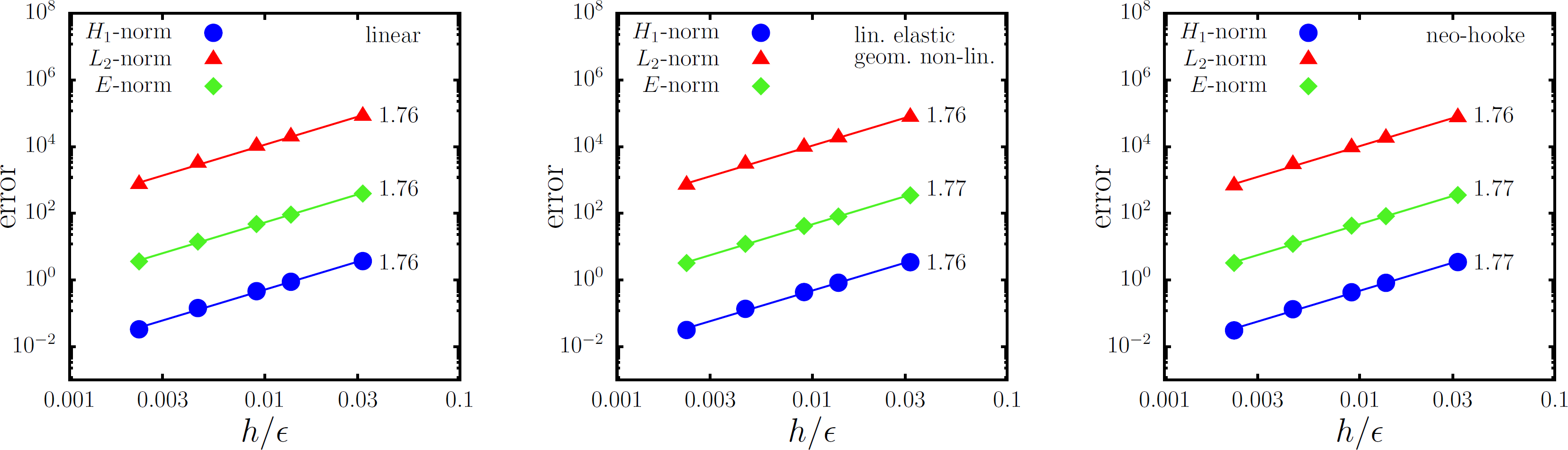} 
	\caption{{\bf Escher's Pegasus}: Micro convergence for (left:) fully linear setting, (center:) linear elasticity along with geometrical nonlinearity, (right:) Neo-Hookean hyperelasticity.
		\label{fig:EscherPegasus2d-Micro-Convergence}}
\end{Figure}

The convergence diagrams in Fig.~\ref{fig:EscherPegasus2d-Micro-Convergence} indicate that the convergence order in all considered norms is virtually insensitive to the nonlinearities introduced by the kinematics and the hyperelastic constitutive law. 
The measured convergence orders are below the full theoretical order of $2q=2$
which is caused by the reduced regularity of the micro problem already visible in the results for the fully linear case.

\subsubsection{Standard Newton versus Alternating Newton} 
\label{subsubsec:SNvAN_pegasus}

In order to calculate the speedup factor of the alternating Newton method related to the standard Newton method, the solution for ndof$_{\text{mac}}=24$ along with ndof$_{\text{mic}}=2\,262\,584$ is calculated by both methods employing 4 loading steps. The corresponding overall performance and computation times are compared in Table~\ref{tab:EscherPegasus2d-SNvAN-LE-micro-mod}. 
  
%--------------------------

\begin{Table}[htbp]
\center
\setlength\tabcolsep{10pt} % default value: 6pt
\begin{tabulary}{\textwidth}{LCCCCC}

%\hline \\[-4.5mm]
%\rowcolor[HTML]{9EBDD5} 
%\multicolumn{5}{l}{\textbf{Standard Newton Method}} \\ 
\hline \\[-4mm]
        
method &        
load step & 
$N_{\text{ite}}^{\text{mac}}$ & 
$t_{\text{ite}}^{\text{mac}}$ &
$u_{\text{max}}$ &
$t_{\text{LS}}$ \\  

\hline \\[-4mm]                       

standard Newton & 1 & 4 & 00:07:27 & 276.8 & 00:27:09 \\
 & 2 & 5 & 00:08:45 & 550.7 & 00:42:31 \\
 & 3 & 5 & 00:08:47 & 818.1 & 00:43:53 \\
 & 4 & 5 & 00:08:45 & 1076.5 & 00:43:52 \\[1mm]

\hline \\[-4mm]

alternating Newton & 1 & 5 & 00:05:08 & 276.8 & 00:25:27 \\
 & 2 & 5 & 00:05:09 & 550.7 & 00:25:45 \\
 & 3 & 5 & 00:05:10 & 818.1 & 00:25:48 \\
 & 4 & 5 & 00:05:09 & 1076.5 & 00:25:52 \\[1mm]         

\hline \\[-4mm]
\multicolumn{6}{r}{total speedup factor = 1.53} \\
\hline \\[-4mm]
\end{tabulary}
\caption{{\bf Escher's Pegasus:} Standard Newton versus alternating Newton with $N_{\text{ite}}^{\text{mac}}$ the number of macro iterations, $t_{\text{ite}}^{\text{mac}}$ the compute time for the first macro Newton iteration including micro Newton iteration(s), $u_{\text{max}}$ the maximum of the normed nodal displacement vectors and $t_{\text{LS}}$ the compute time for a complete load step. Compute time has the format (HH:mm:ss) for hours, minutes, and seconds.} 
\label{tab:EscherPegasus2d-SNvAN-LE-micro-mod}
\end{Table}
  
Table~\ref{tab:EscherPegasus2d-SNvAN-LE-micro-mod} indicates that the alternating Newton converges to the correct solution in terms of the maximum nodal displacement. Remarkably, in only one out of four loading steps the alternating Newton requires one additional macro Newton iteration for convergence. The overall speedup of $1.53$ indicates the improved efficiency of the alternating Newton.

\begin{Table}[htbp]
\center
\begin{tabular}{ccrr}
\hline \\[-4mm] 

load step number & macro ite & \multicolumn{2}{c}{macro residual}    \\
          &           & standard Newton  & alternating Newton \\ 
\hline \\[-4mm]                       
\small
4 (out of 4) & 1 & 1.0000e+00 & 9.9999e-01 \\
             & 2 & 1.4152e+01 & 1.4114e+01 \\
             & 3 & 6.3705e-02 & 4.7410e-02 \\
             & 4 & 4.5369e-04 & 5.8777e-04 \\
             & 5 & 1.8262e-08 & 6.1250e-06 \\[1mm]
\hline
\small
1 (out of 1) & 1 & 1.0000e+00 & 1.0000e+00 \\
             & 2 & 6.6287e+01 & 6.6301e+01 \\
             & 3 & 4.0499e+00 & 3.1304e+00 \\
             & 4 & 2.3008e-02 & 1.1509e-02 \\
             & 5 & 7.5719e-05 & 1.4101e-04 \\
		     & 6 & 5.0931e-09 & 1.5738e-06 \\[1mm]
\hline \\[-4mm]
\end{tabular}
\caption{{\bf Escher's Pegasus:} Comparison of the macro residuals for standard Newton versus alternating Newton in load step 4 (out of 4) and load step 1 (out of 1).} 
\label{tab:EscherPegasus2d-SNvAN-convergence-LE-micro}
\end{Table}
  
Table~\ref{tab:EscherPegasus2d-SNvAN-convergence-LE-micro} compares the convergence of the macro Newton for the two different solution strategies in terms of their macro residuals. {\color{black} For that aim, two cases are considered; in the first case the total load is decomposed in four increments of equal magnitude, in the second case in one single step.} Independent of the magnitude of the load step the number of required iterations is the same for both methods, although the alternating Newton exhibits a minor deviation from quadratic convergence in its last iteration. 

\begin{Table}[htbp]
\center
\renewcommand{\arraystretch}{1.1}
\begin{tabular}{rrcccc}
\hline \\[-5mm]         
macro       &  micro       &  \multicolumn{4}{c}{speedup factor}   \\
  ndof      &  ndof        &  \multicolumn{2}{c}{linear elastic}  &  \multicolumn{2}{c}{neo-Hookean}   \\
	          &              &  \small{$N_{\text{LS}}=$4} & \small{$N_{\text{LS}}=$1} & \small{$N_{\text{LS}}=$4} & \small{$N_{\text{LS}}=$1} \\ 
\hline \\[-5mm]                       

24   &    1\,946     &  1.50 & 1.38 & 1.50 & 1.27 \\
     &    11\,942    &  1.40 & 1.56 & 1.31 & 1.28 \\
     &    25\,820    &  1.44 & 1.48 & 1.39 & 1.29 \\
     &    101\,918   &  1.44 & 1.47 & 1.38 & 1.31 \\
     &    384\,930   &  1.48 & 1.48 & 1.39 & 1.30 \\
     &    2\,262\,584  &  1.53 & 1.56 & 1.45 & 1.40 \\
\hline \\[-5mm] 

1\,122      &    1946 &  1.64 & 1.63 &  1.35 & 2.57    \\
4\,242      &         &  1.66 & 1.58 &  1.31 & 1.80 \\
16\,482     &         &  1.59 & 1.59 &  1.31 & 1.83 \\
64\,962     &         &  1.56 & 1.40 &  1.33 & 1.76 \\
257\,922    &         &  1.57 & \color{black} 1.44 &  1.33 & 1.70  \\
4\,103\,682   &         &  1.55 & \color{black} 1.49 &  1.31 & 1.63
\\ 
\hline
\end{tabular} 
\caption{{\bf Escher's Pegasus:} speedup factors of alternating Newton compared to standard Newton depending on discretizations and the type of the employed elastic constitutive law, where $N_{\text{LS}}$ is the total number of load steps.
\label{tab:Escher-pegasus-2d-speedup-MicVariation}}
\end{Table}
\begin{Figure}[htbp]
	\begin{minipage}{16.5cm}  
		\centering                                                 
		\includegraphics[height=5.5cm, angle=0, clip=]{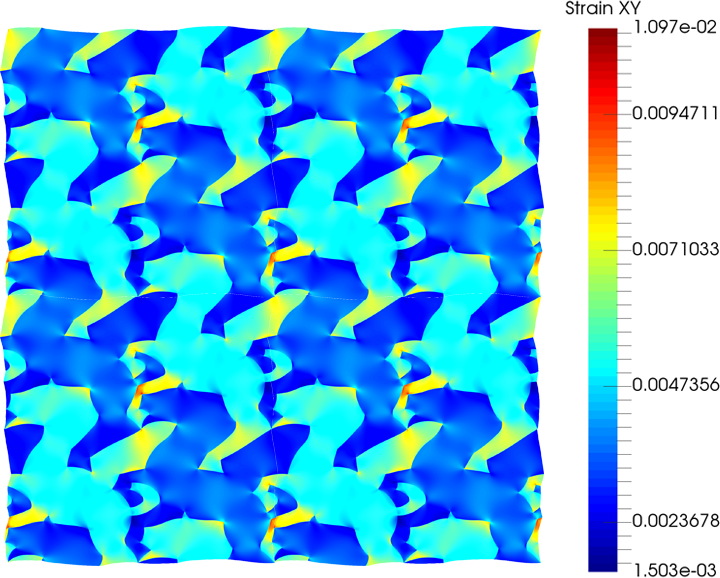}
		\hspace*{4mm}
		\includegraphics[height=5.5cm, angle=0, clip=]{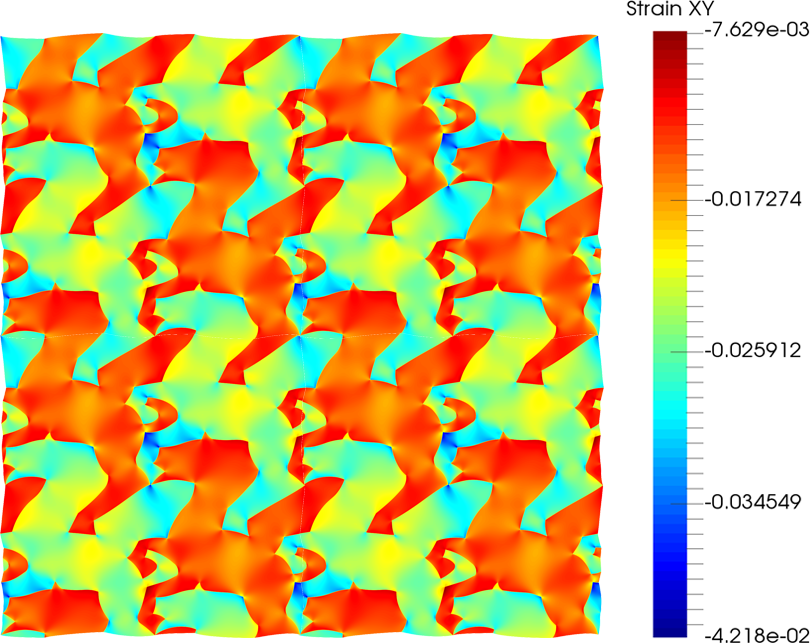}
	\end{minipage}
	\caption{{\bf Escher's Pegasus}: contour plot of shear strain component $E_{xy}$ (left) at $x=3943.37 \text{mm}$, $y=211.32$~mm and (right) at $x=1056.62$~mm, $y=788.67$~mm.
		\label{fig:Escher-pegasus-colored-v1}}
\end{Figure} 
Table~\ref{tab:Escher-pegasus-2d-speedup-MicVariation} lists the speedup factors of the novel solution strategy for different discretizations; they vary from 1.27 to almost 2.6. Increasing the micro discretization while keeping the macro discretization fixed results in a modest increase of the speedup factor. Vice versa, increasing the macro discretization while keeping the micro discretization fixed results in a modest decrease of the speedup factors. The two trends hold --apart from some outliers-- for linear elasticity as well as for neo-Hookean hyperelasticity and they are insensitive to the load step size likewise. The speedup factor for linear elasticity is generally higher than for the neo-Hookean case, the exception is the case where ndof$_{\text{mac}}$ is in the order of or larger than ndof$_{\text{mic}}$; here the speedup factor for the hyperelastic case is higher, if the load is applied in one single step instead of four steps. 
 
Figure~\ref{fig:Escher-pegasus-colored-v1} shows the distribution of shear strain $2E_{xy}$ for a $2 \times 2$ RVE-array of the Pegasus microstructure with periodic fluctuations at two selected points in the beam. 
  
%-----------------------------
\subsection{Escher's Swan} 
\label{subsec:Escher_swans}
%-----------------------------

\begin{Figure}[htbp]
	\begin{minipage}{16.0cm}  
		\centering 
            \includegraphics[height=5.0cm, angle=0, clip=]{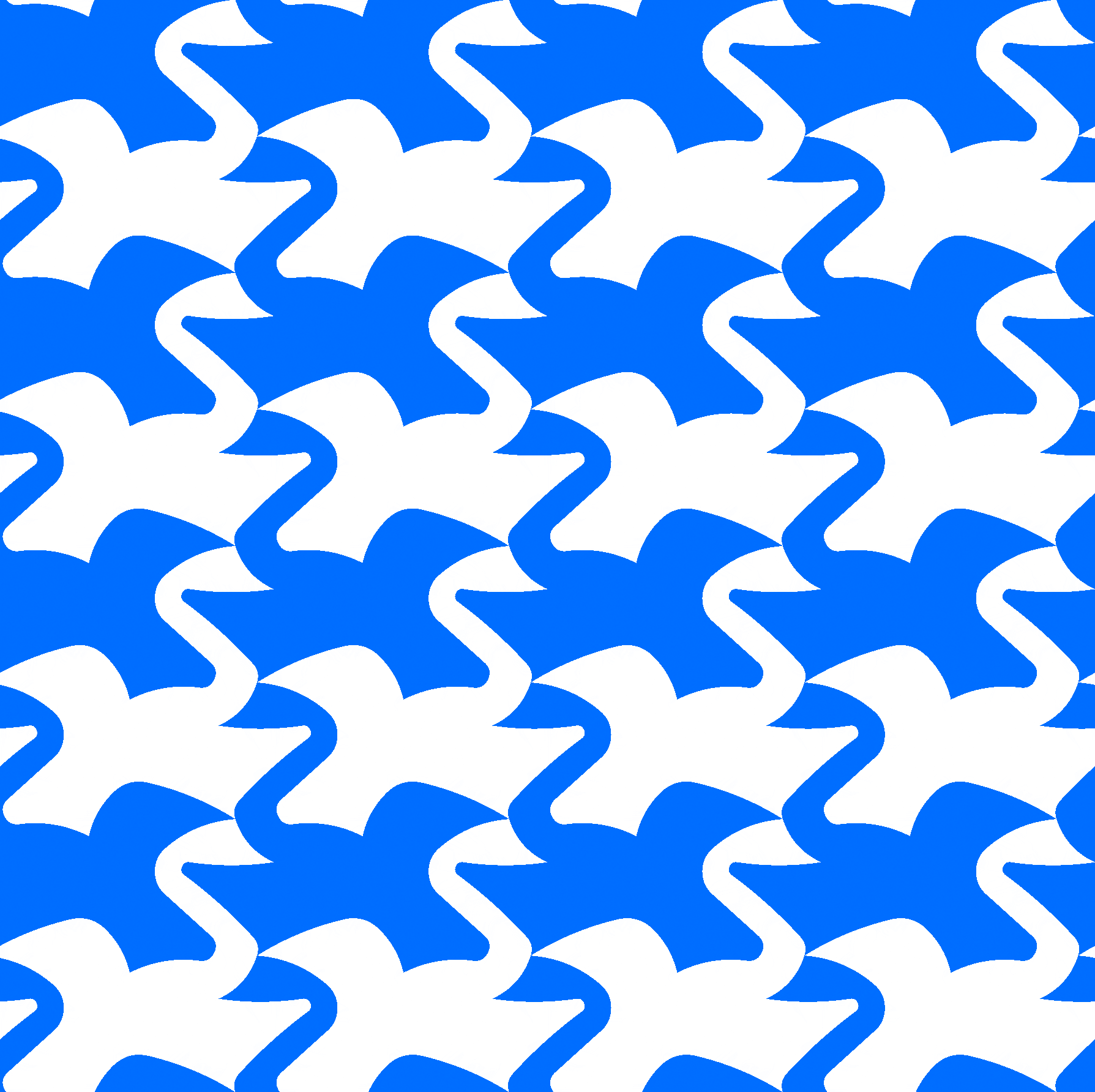} 
		\hspace*{8mm}
	    	\includegraphics[height=5.0cm, angle=0, clip=]{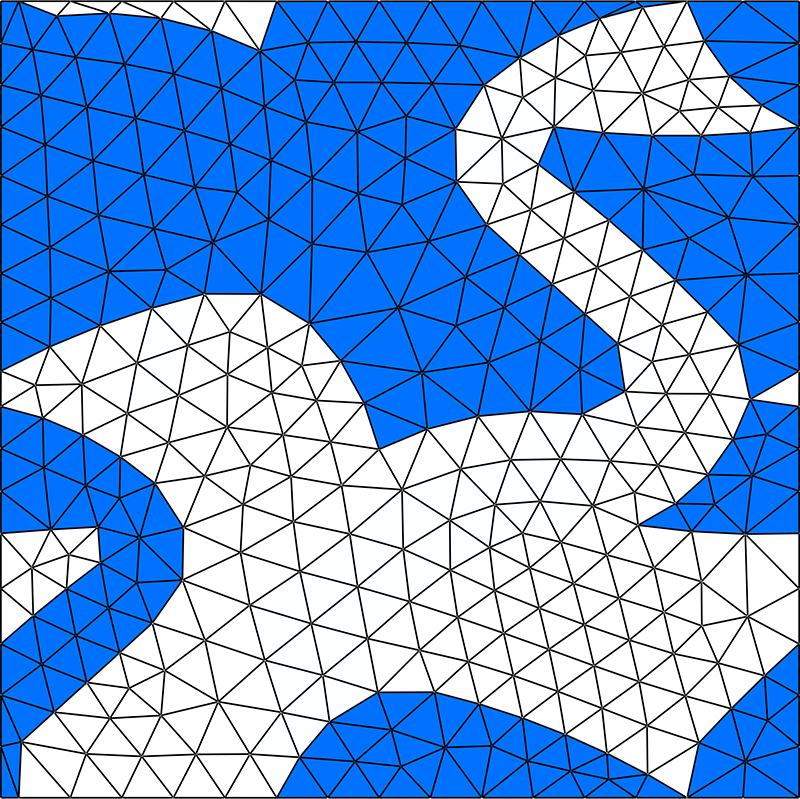}
	\end{minipage}
	\caption{{\bf Escher's Swan:}  (left) Periodic tessellation, (right) discretization of a unit cell.
		\label{fig:Escher-Swans}}
\end{Figure}

In this example a two-phase microstructure is used that follows Escher's ''Swan N$^{\circ}$96'' tessellation as displayed in the left of Fig.~\ref{fig:Escher-Swans}; each phase follows from each other by reflection and translation. Figure~\ref{fig:Escher-Swans} (right) shows a unit cell along with a discretization using linear triangles, hence $q=1$. 
In this example phase 1 is assigned to the swan to the left in blue color, phase 2 is assigned to the swan to the right in white color. The material parameters for linear elastic and neo-Hookean material laws are given in Tab.~\ref{tab:Material-param-pegasus} of Sec.~\ref{subsec:Escher_pegasus}.

The macro scale BVP is again the cantilever beam with unaltered dimensions $l, b, t$ which is subject at $x=l$ to a displacement control where the total displacement of $1100$~mm is applied in four loading steps.
%------------------------------------------

\subsubsection{Micro convergence for linear elasticity and neo-Hookean hyperelasticity}

\begin{Figure}[htbp]
	\centering
	\includegraphics[height=4.5cm, angle=0, clip=]{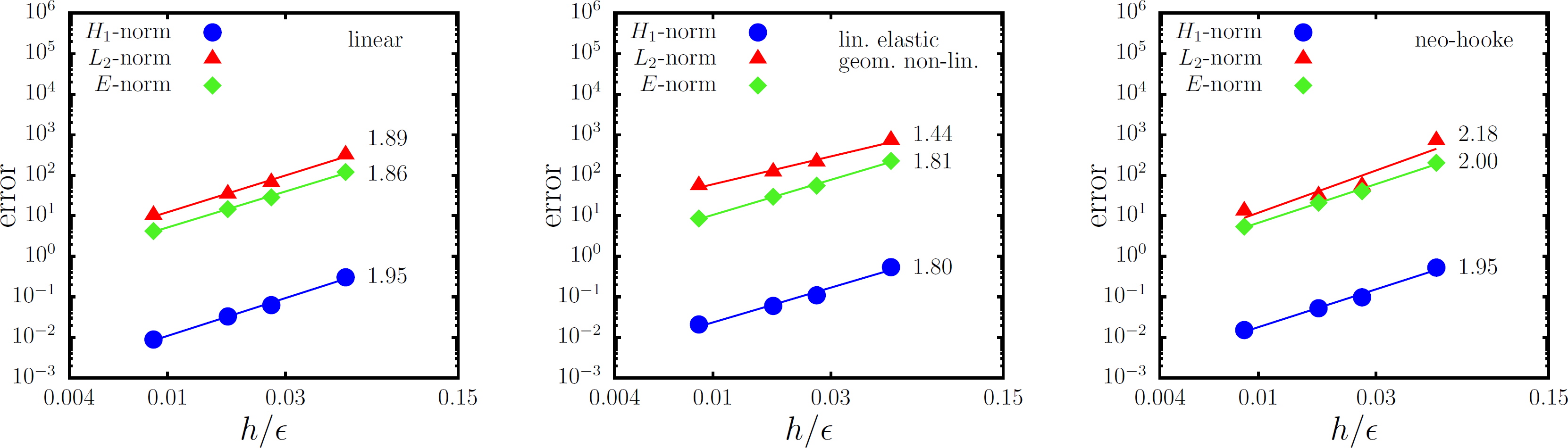} 
	\caption{{\bf Escher's Swan:} Micro convergence for (left:) the fully linear case, (center:) for linear elasticity with geometrical nonlinearity, and (right:) for Neo-Hookean hyperelasticity.
		\label{fig:EscherSwan2d-Micro-Convergence}}
\end{Figure}
The reference solution for the micro convergence analysis is obtained for ndof=24 on the macro domain (5$\times$1 macro elements) and ndof=574\,922 on micro domains with $\epsilon=11.4$~mm.

The convergence diagrams for linear elasticity and neo-Hookean material laws are displayed in Fig.~\ref{fig:EscherSwan2d-Micro-Convergence}. Here, the fully linear case achieves almost the full theoretical order in all norms, $2q=2$. For the combination of linear elasticity and  geometrical nonlinearity, the convergence order is throughout reduced compared to the linear case, most pronounced in the $L_2$-norm.
Remarkably, the convergence order for the fully nonlinear case of neo-Hookean hyperelasticity exhibits excellent convergence rates, in the $L_2$-norm even above the nominal order 2.0.
  
\subsubsection{Standard Newton versus alternating Newton}
\label{subsubsec:SNvAN_study2}

\begin{Table}[htbp]
	\center
	\setlength\tabcolsep{10pt} % default value: 6pt
	\begin{tabulary}{\textwidth}{LCCCCC}

		\hline \\[-4mm]
		
		method &        
		load step & 
		$N_{\text{ite}}^{\text{mac}}$ & 
		$t_{\text{ite}}^{\text{mac}}$ &
		$u_{\text{max}}$ &
		$t_{\text{LS}}$ \\  
		
		\hline \\[-4mm]                       
		
		standard Newton & 1 & 4 & 00:01:54 &  279.3 & 00:07:03 \\
		& 2 & 4 & 00:02:12 &  561.7 & 00:08:32 \\
		& 3 & 4 & 00:02:12 &  847.3 & 00:08:49 \\
		& 4 & 4 & 00:02:16 & 1136.6 & 00:08:59 \\[1mm]
		
		\hline \\[-4mm]
		
		alternating Newton & 1 & 4 & 00:01:22 &  279.3 & 00:05:27 \\
		& 2 & 5 & 00:01:23 &  561.7 & 00:06:54 \\
		& 3 & 5 & 00:01:23 &  847.3 & 00:06:54 \\
		& 4 & 5 & 00:01:22 & 1136.6 & 00:06:53 \\[1mm]         
		
		\hline \\[-4mm]
		\multicolumn{6}{r}{total speedup factor = 1.28} \\
		\hline \\[-4mm]
	\end{tabulary}
	\caption{{\bf Escher's Swan:} Standard Newton versus alternating Newton where $N_{\text{ite}}^{\text{mac}}$ is the number of macro iterations, $t_{\text{ite}}^{\text{mac}}$ the compute time for the first macro Newton iteration including micro Newton iteration(s), $u_{\text{max}}$ is the maximum of the normed nodal displacement vectors and $t_{\text{LS}}$ is the compute time for a load step.} 
	\label{tab:EscherSwan2d-SNvAN-LE-micro-mod}
\end{Table}

\begin{Table}[htbp]
	\center
	\renewcommand{\arraystretch}{1.1}
	\begin{tabular}{rrcccc}
		\hline \\[-5mm]         
		macro       &  micro       &  \multicolumn{2}{c}{speedup factor}   \\
		ndof      &  ndof        &  linear elastic  &  neo-Hookean   \\
		\hline \\[-5mm]                       
		
		24   &    674     &  1.25 & 1.20  \\
		&    3\,248    &  1.68 & 1.05  \\
		&    7\,222    &  1.26 & 1.11  \\
		&    27\,444   &  1.25 & 1.20  \\
		&    74\,486   &  1.21 & 1.18  \\
		&    574\,922  &  1.28 & 1.23  \\
		\hline \\[-5mm] 
		
		1\,122      &    674  &  2.13 & 1.50 \\
		4\,242      &         &  1.53 & 1.44 \\
		16\,482     &         &  1.56 & 1.26 \\
		64\,962     &         &  1.56 & 1.37 \\
		257\,922    &         &  1.76 & 1.44 \\
		4\,103\,682   &         &  \color{black} 1.63 & 1.49 
		\\ 
		\hline
	\end{tabular}
	\caption{{\bf Escher's Swan:} speedup factors of alternating Newton compared to standard Newton for different discretizations and types of elastic constitutive laws.
		\label{tab:Escher-swan-2d-speedup-MicVariation}}
\end{Table}
  
The overall performance and computational time of the reference solution calculated by standard and alternating Newton method is compared in Tab.~\ref{tab:EscherSwan2d-SNvAN-LE-micro-mod}. 
Maximum nodal displacements are in agreement for both methods, the speedup factor is 1.28. In the present example the alternating Newton method requires in 3 out of 4 loading steps 5 macro Newton iterations and therefore one more than for the standard, nested Newton.
 
A richer picture on the speedup factors of the alternating Newton method for different discretizations is given in Tab.~\ref{tab:Escher-swan-2d-speedup-MicVariation}. Among the four cases considered therein the case of linear elasticity along with geometrical nonlinearity results for ndof$_{\text{mac}}$ $\gg$ ndof$_{\text{mic}}$ in the highest speedup factors, throughout above 1.5.   

\begin{Figure}[htbp]
	\begin{minipage}{16.0cm}  
		\centering
		\includegraphics[height=6.0cm, angle=0, clip=]{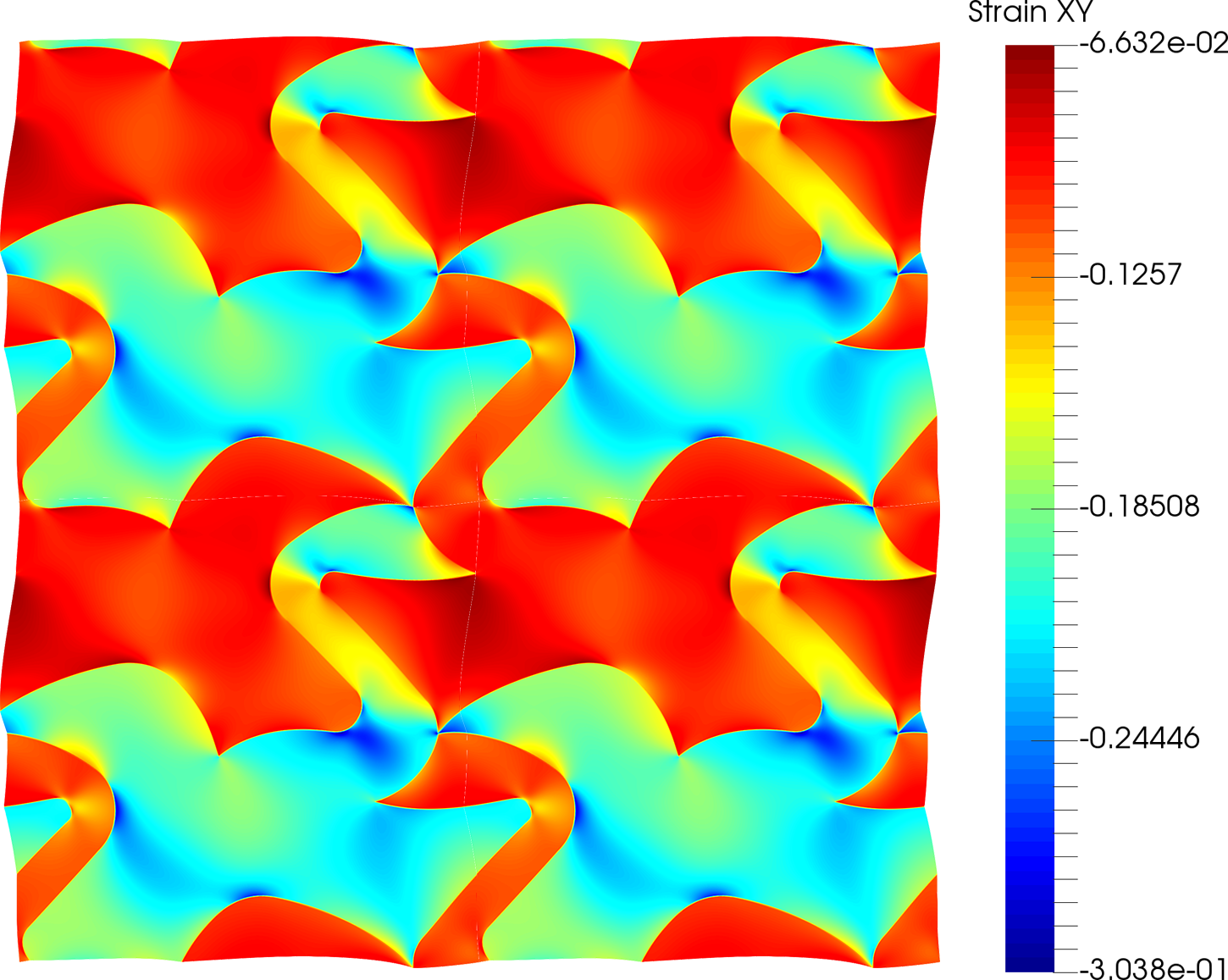}
	    \includegraphics[height=6.0cm, angle=0, clip=]{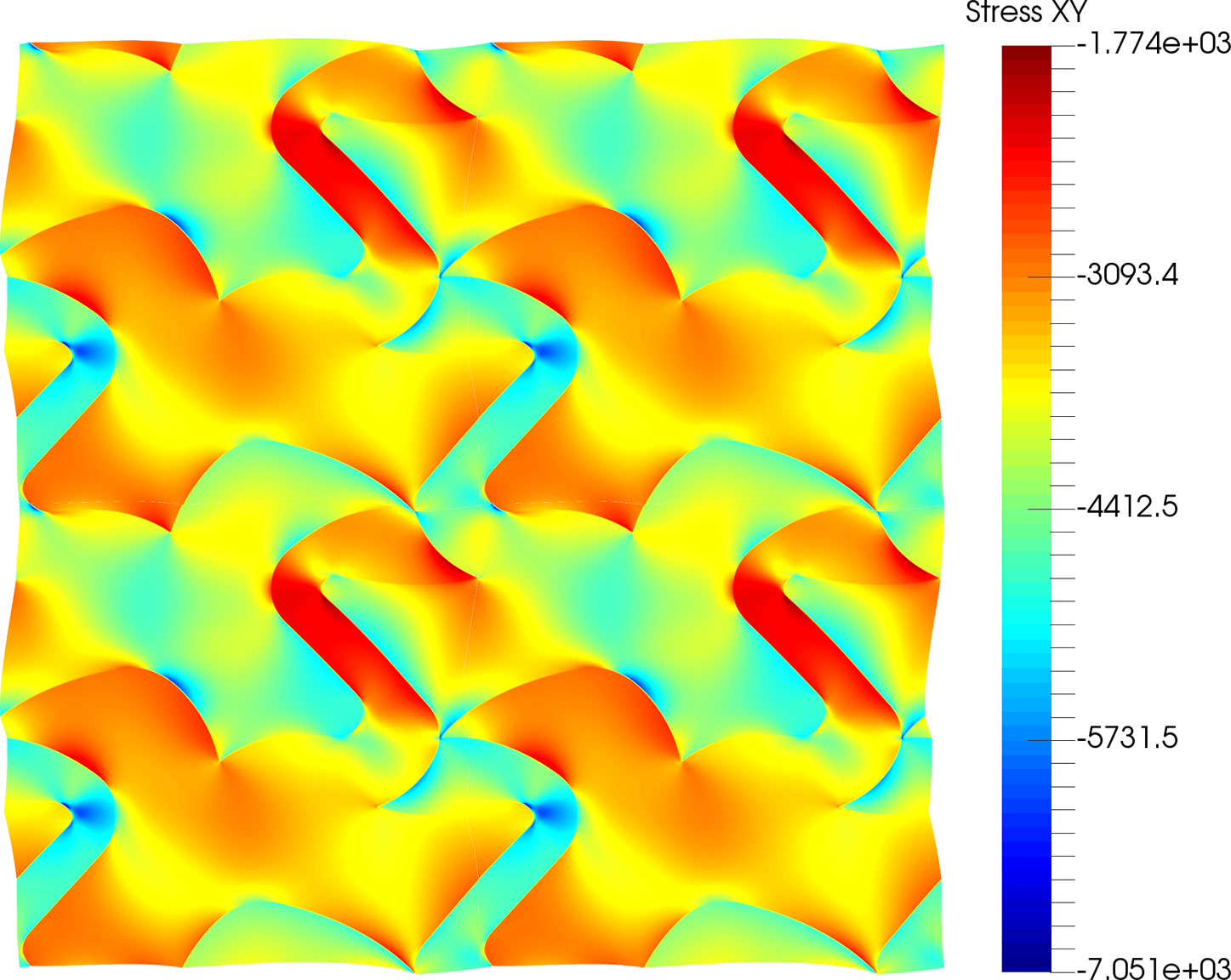}
	\end{minipage}
	\caption{{\bf Escher's Swan:} at $x=1056.62~\text{mm}$, $y=788.67~\text{mm}$ the contour plot (left) of shear strain $2E_{xy}$ and (right) of von-Mises stress $S_{xy}$.
		\label{fig:Escher-swans-GP3-2x2}}
\end{Figure}

Figure~\ref{fig:Escher-swans-GP3-2x2} (left) underpins that  phase 1 (swan to the left) is the stiffer one due to its smaller shear strains, the contour plot in the right displays the corresponding shear stress distribution.
 
%-----------------------------
\subsection{DFG-Heisenberg} 
\label{subsec:DFG_Heisenberg}
%-----------------------------
 
\begin{Figure}[htbp]
   \begin{minipage}{16.0cm}  
   %\centering
         \hspace*{10mm}
         \includegraphics[height=6.0cm, angle=0, clip=]{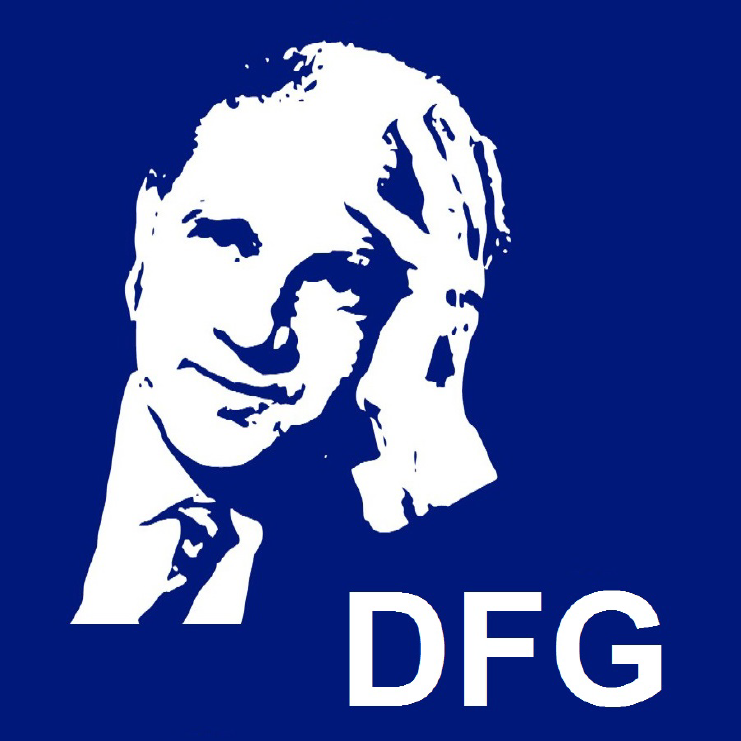}  \hspace*{4mm}
         \includegraphics[height=6.30cm, angle=0,  clip=]{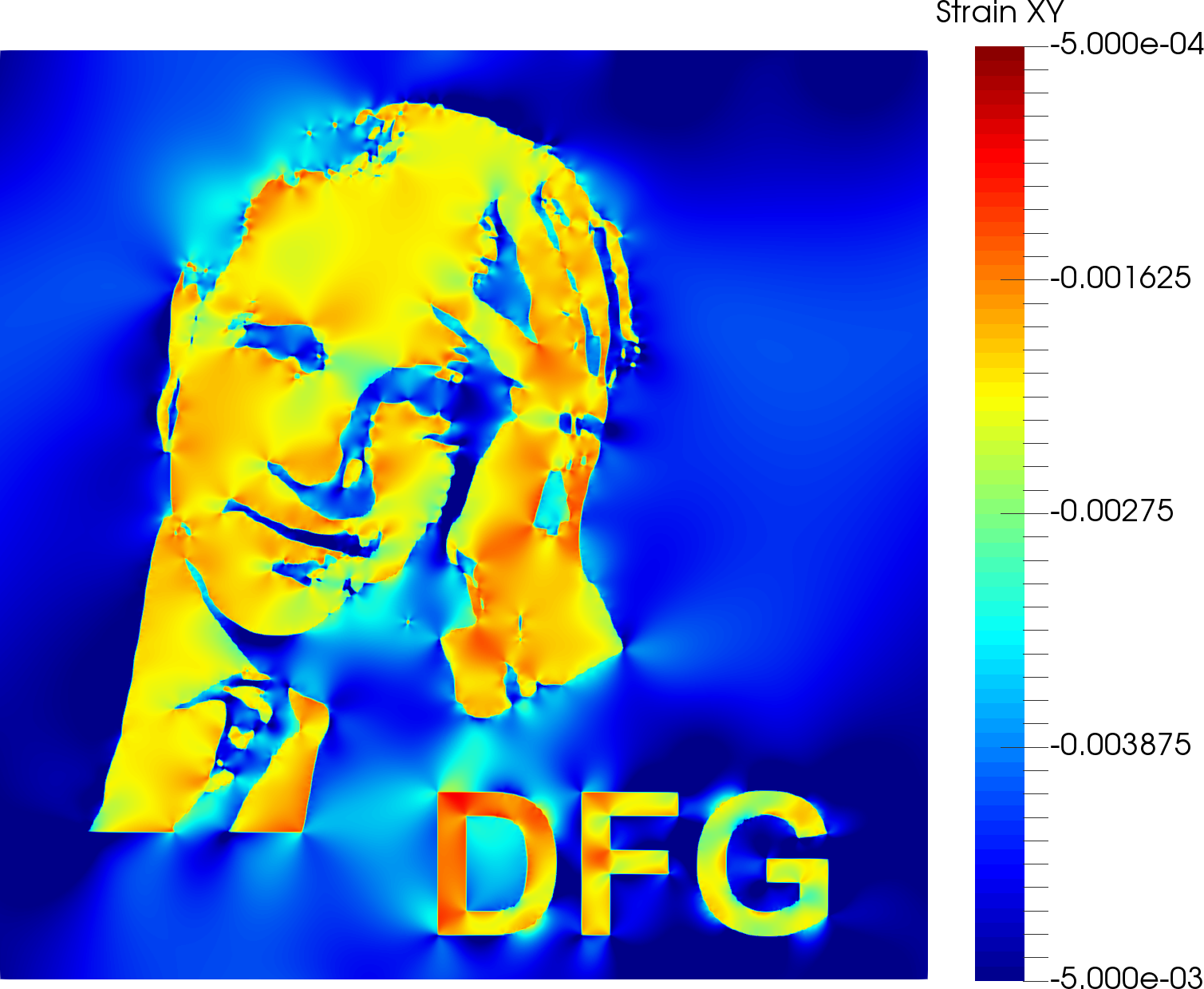}
   \end{minipage}
\caption{{\bf DFG-Heisenberg}: (Left) two-phase distribution, (right) contour plot of shear strain $2E_{xy}$ at position $x=666.67$~mm, $y=333.34$~mm. 
\label{fig:DFG_Heisenberg}}
\end{Figure}
 
As an example of pixel-based microstructure representation we consider the two-phase portrait of the physicist Werner Heisenberg\footnote{in the format that serves as the program logo of the German Research Foundation (DFG).} as shown in the left of Fig.~\ref{fig:DFG_Heisenberg}. 
 
The discretization of this image is carried out using pixel to mesh conversion, phases are defined by conversion into binaries according to the pixels' rgb color code. The blue matrix is the softer phase 1, the white inclusions are the stiffer phase 2 as indicated by the reduced strain in the right of  Fig.~\ref{fig:DFG_Heisenberg}. The elasticity parameters of each phase are listed in Tab.~\ref{tab:Material-param-pegasus}. 
 
The macro scale BVP is a cantilever beam with length $l=1000$~mm, $b=1000$~mm, $t=100$~mm, which is loaded at $x=l$ by a line load of magnitude $100$~N/mm in negative $y$-direction. PBCs are applied to the micro domains of length $\epsilon=1$~mm.

\subsubsection{Micro convergence for linear elasticity and neo-Hookean hyperelasticity}

\begin{Figure}[htbp]
	\centering
	\includegraphics[height=4.5cm, angle=0, clip=]{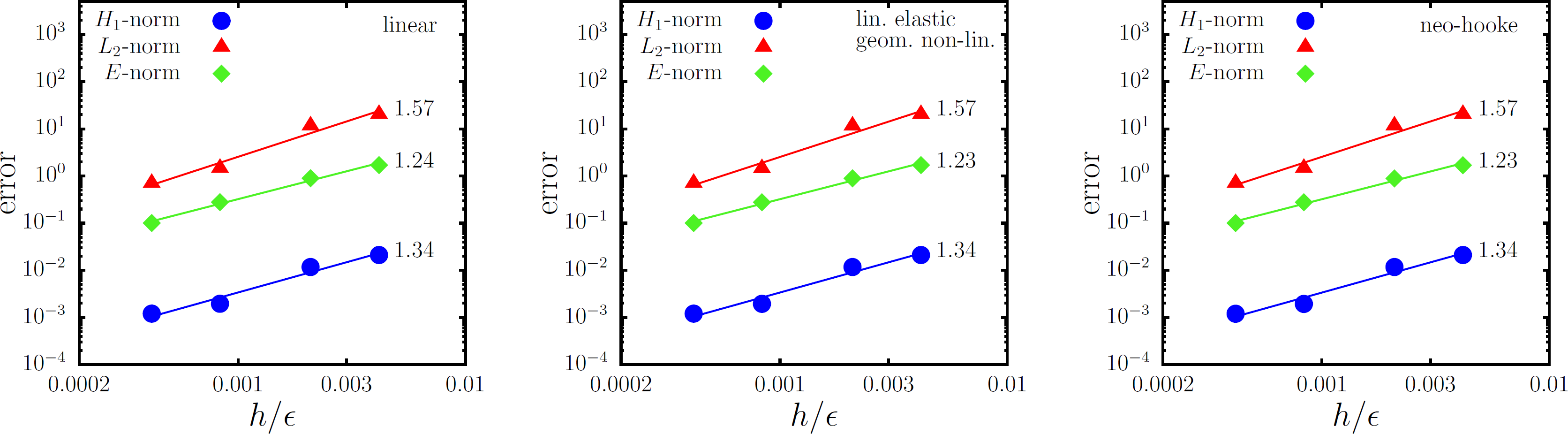} 
	\caption{{\bf DFG-Heisenberg:} Micro convergence for (left:) the fully linear case, (center:) for linear elasticity with geometrical nonlinearity, and (right:) for neo-Hookean hyperelasticity.
		\label{fig:Heisenberg2d-Micro-Convergence}}
\end{Figure}
The reference solution for the convergence analysis is obtained by a discretization of ndof= 46\,099\,202 on the micro domain employing a characteristic element width of $h=1/4800$~mm with $q=1$. The discretization directly follows from the pixels of the original image. Since we are mostly interested in the micro convergence and differences between the cases having different levels of nonlinearity, a simple driver on the macro scale in terms of one triangular macro element with $p=1$ is enough.
  
The measured convergence orders in the diagrams of Fig.~\ref{fig:Heisenberg2d-Micro-Convergence} are considerably reduced compared to the theoretical order of $2q=2$; it is most pronounced in the energy-norm and least pronounced in the $L_2$-norm. The reduced order is caused by the stiff inclusion phase pinching into the softer matrix phase, compare the shear strain distribution in Fig.~\ref{fig:DFG_Heisenberg}. Due to almost identical results for the fully linear case with the nonlinear cases it can be concluded that the nonlinearities by geometry and material have no impact on convergence. 
 
\subsubsection{Standard Newton versus alternating Newton}
\label{subsubsec:SNvAN_heisenberg}
 
The speedup factors of the novel scheme compared to the standard solution procedure are given in  Tab.~\ref{tab:DFG-heisenberg-2d-speedup-MicVariation}.

\begin{Table}[htbp]
	\center
	\renewcommand{\arraystretch}{1.1}
	\begin{tabular}{rrcccc}
		\hline \\[-5mm]         
		macro       &  micro       &  \multicolumn{2}{c}{speedup factor}   \\
		ndof      &  ndof        &  linear elastic  &  neo-Hookean   \\
		\hline \\[-5mm]                       
		
		8   &    116162     &  1.41 & 1.41  \\
		&    462\,722     &  1.37 & 1.32  \\
		&    1\,847\,042    &  1.51 & 1.41  \\
		&    2\,884\,802    &  1.58 & 1.49  \\
		&    11\,529\,602   &  1.82 & 1.58  \\
 		&    46\,099\,202   &  1.73 & 1.73  \\
		\hline
		882     &    116\,162     &  1.21 & 1.38  \\
		3\,362    &               &  1.22 & 1.41  \\
		13\,122   &               &  1.22 & 1.41  \\
		51\,842   &               &  1.28 & 1.41  \\
 		206\,082  &               &  1.27 & 1.25  \\
		\hline
	\end{tabular}
	\caption{{\bf DFG-Heisenberg:} speedup factors of alternating Newton compared to standard Newton depending on discretizations and type of elastic constitutive law. \label{tab:DFG-heisenberg-2d-speedup-MicVariation}}
\end{Table}

The speedup is throughout in between 1.2 and 1.8, such that the novel solution strategy pays off in either case. For a fixed macro discretization the speedup factor is the larger, the finer the micro discretization, which is independent of the elasticity law and true apart from some outliers. For the case of a fixed micro discretization the speedup-factor is relatively insensitive to the macro discretization. These trends are consistent with the previous example of Escher's Pegasus.
 
%---------------------------------------------------------------------------------------------------------
\subsection{Nanoporous gold}
\label{subsec:Nanoporous-gold}
%---------------------------------------------------------------------------------------------------------

\begin{Figure}[htbp]
	\begin{minipage}{16.0cm}  
		\centering
		\includegraphics[height=6.2cm, angle=0, clip=]{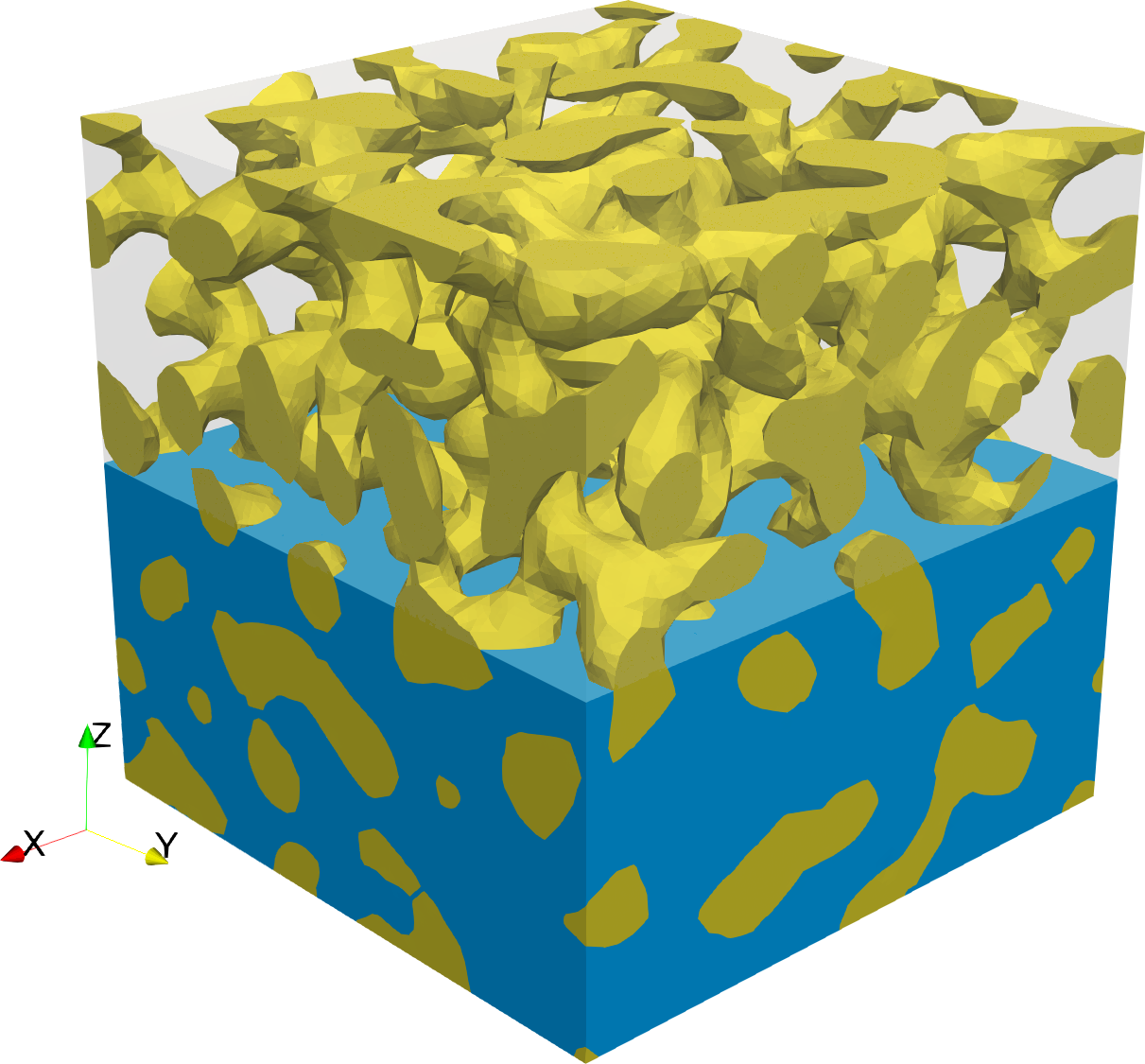}  \hspace*{4mm}
		\includegraphics[height=6.2cm, angle=0, clip=]{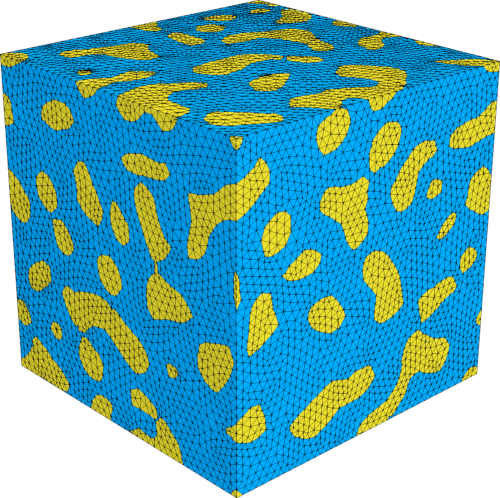}  \hspace*{10mm}
	\end{minipage}
	\caption{{\bf Nanoporous gold composite}: (Left) gold phase distribution, (right) discretization with ndof=1\,807\,389.
		\label{fig:Nanoporous_gold}}
\end{Figure}
  
Nanoporous gold exhibits extraordinary properties such as high catalytic activity \cite{Fujita-etal-2012}, surface-stress induced macroscopic bending of cantilevers  \cite{Kramer-etal2004}, and it is a promising component for gold-polymer nanocomposites \cite{GriffithsBargmannReddy2017}.

Here we consider a gold-polymer nanocomposites with the nanoporous gold phase embedded in a softer second phase. Linear Dirichlet boundary conditions are applied on the RVE. Instead of using true material parameters for the phases, we stick to the unified material parameters valid throughout this work. The Young's moduli of the gold phase are set to $E_i=100\,000$ $\text{N/mm}^2$ and the matrix phase $E_m=40\,000$ $\text{N/mm}^2$, $\nu=0.2$ for both phases, see Tab. \ref{tab:Material-param-pegasus}. The volume ratio of the gold phase is $V_i/V_{\text{tot}}=0.36$.   

The macro BVP is a cantilever beam with length $l=5000~\text{mm}$, height $b=1000~\text{mm}$, and thickness $t=1000~\text{mm}$. At its free end the cantilever beam is subject to displacement control with a total displacement of $200~\text{mm}$ applied to all nodes at $x=l$ in four loading steps. 

\subsubsection{Micro convergence for linear elasticity and neo-Hookean hyperelasticity}
The reference solution for the convergence study is obtained by ndof=72 on the macro domain (5$\times$1$\times$1 hexahedral macro elements) and ndof=5\,484\,291 on the micro domain with $\epsilon=4.03$~mm.

\begin{Figure}[htbp]
	\centering
	\includegraphics[height=4.5cm, angle=0, clip=]{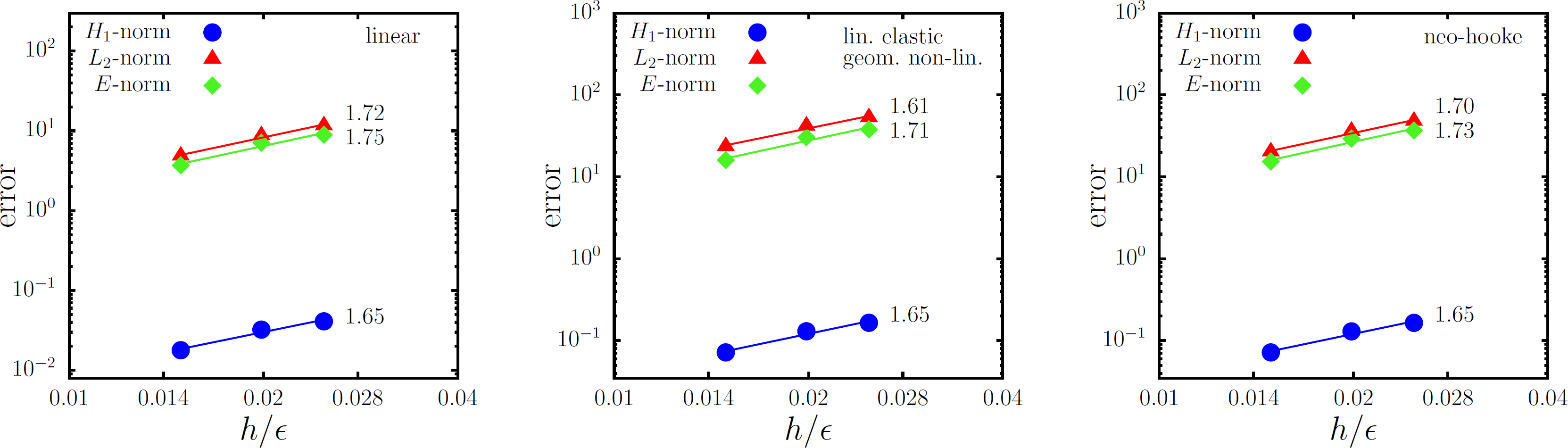} 
	\caption{{\bf Nanoporous gold composite}: Micro convergence for (left:) the fully linear case, (center:) linear elasticity with geometrical nonlinearity, and (right:) the case of Neo-Hookean hyperelasticity.
		\label{fig:Nanogold-Micro-Convergence}}
\end{Figure}
  
The convergence diagrams are displayed in Fig.~\ref{fig:Nanogold-Micro-Convergence}. The convergence order for the case of geometrical nonlinearity as well as for additional material nonlinearity is virtually the same as for the fully linear case, an exception is the case of geometrical nonlinearity along with linear elasticity in the $L_2$-norm.
 
\subsubsection{Standard Newton versus alternating Newton}
\label{subsubsec:SNvAN_NanoporousGold}

\begin{Table}[htbp]
\center
\renewcommand{\arraystretch}{1.1}
\begin{tabular}{rrcccc}
\hline \\[-5mm]         
macro       &  micro       &  \multicolumn{2}{c}{speedup factor}   \\
  ndof      &  ndof        &  linear elastic  &  neo-Hookean   \\
\hline \\[-5mm]                       

 72   &    262\,308       &  1.75           & 1.69  \\
      &    353\,904       &  1.77           & 1.71  \\
      &    690\,162       &  1.98           & 1.82  \\
      &    1\,518\,654    &  1.88           & 1.86  \\
      &    1\,807\,389    &  1.86           & 1.86  \\
\hline
 297  &    262\,308    &  1.85           & 1.85  \\
 1575 &                &  1.90           & 1.89  \\
 4557 &                &  1.78           & 1.81  \\
\hline
\end{tabular}
\caption{{\bf Nanoporous gold composite}: speedup factors of the alternating Newton compared to the standard Newton depending on discretizations and the type of elasticity law. 
\label{tab:Nanoporous-gold-RVE-3d-speedup-MicVariation}}
\end{Table}

The speedup factors in Tab.~\ref{tab:Nanoporous-gold-RVE-3d-speedup-MicVariation} range from 1.7 to almost 2 and are higher than in the previously considered homogenization problems in 2d. There is no clear trend in the speedup factors for finer meshes except of for the case of neo-Hookean hyperelasticity, where the speedup is the larger, the finer the micro mesh. For the other cases of various discretizations (combinations of micro to macro) the speedup is around 1.8 with relatively minor variance. 

\begin{Figure}[htbp]
	\begin{minipage}{16.5cm}  
		\centering
		\includegraphics[height=5.8cm, angle=0, clip=]{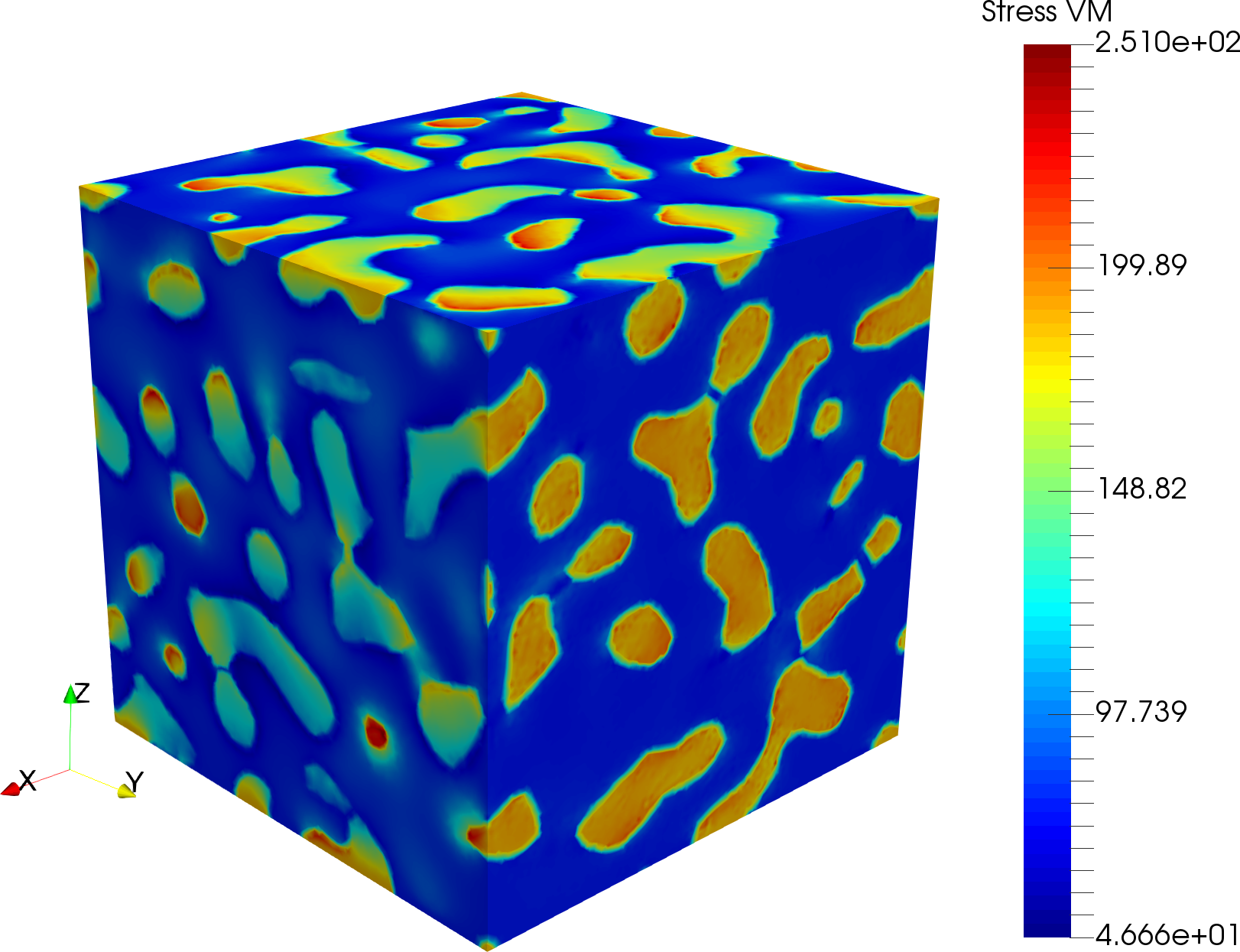}  \hspace*{2mm}
		\includegraphics[height=5.8cm, angle=0, clip=]{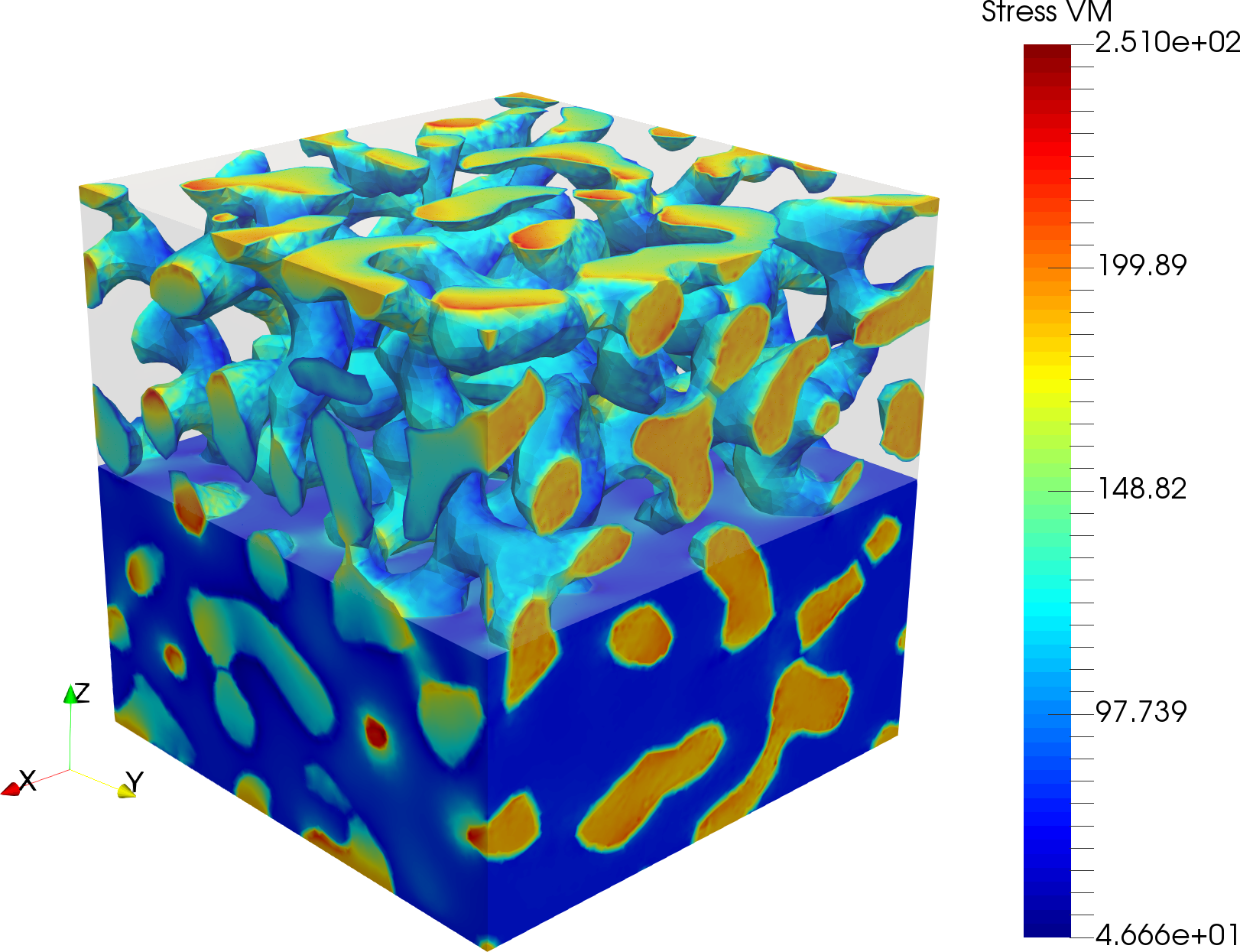}  \hspace*{10mm}
	\end{minipage}
	\caption{{\bf Nanoporous gold composite}: von-Mises stress distribution in (N/mm$^2$) {\color{black} at  $x=y=z=605.6624$~mm}.
		\label{fig:Nanoporous_gold_VM_Stress}}
\end{Figure}

%% file: main.bbl
\begin{thebibliography}{10}

\bibitem{Assyr2006}
Assyr Abdulle.
\newblock Analysis of the heterogeneous multiscale \textsc{FEM} for problems in
  elasticity.
\newblock {\em Mathematical Models and Methods in Applied Sciences},
  16(04):615--635, 2006.

\bibitem{EidelFischer2018}
Bernhard Eidel and Andreas Fischer.
\newblock The heterogeneous multiscale finite element method for the
  homogenization of linear elastic solids and a comparison with the
  \textsc{FE}2 method.
\newblock {\em Computer Methods in Applied Mechanics and Engineering},
  329:332--368, 2018.

\bibitem{MichelMoulinecSuquet1999}
J.~C. Michel, H.~Moulinec, and P.~Suquet.
\newblock Effective properties of composite materials with periodic
  microstructure: a computational approach.
\newblock {\em Computer Methods in Applied Mechanics and Engineering},
  172(1-4):109--143, 1999.

\bibitem{Miehe-etal-1999a}
Christian Miehe, J{\"o}rg Schr{\"o}der, and Jan Schotte.
\newblock Computational homogenization analysis in finite plasticity simulation
  of texture development in polycrystalline materials.
\newblock {\em Computer Methods in Applied Mechanics and Engineering},
  171(3-4):387--418, 1999.

\bibitem{Miehe-etal-1999b}
Christian Miehe, Jan Schotte, and J{\"o}rg Schr{\"o}der.
\newblock Computational micro--macro transitions and overall moduli in the
  analysis of polycrystals at large strains.
\newblock {\em Computational Materials Science}, 16(1-4):372--382, 1999.

\bibitem{FeyelChaboche2000}
Fr{\'e}d{\'e}ric Feyel and Jean-Louis Chaboche.
\newblock \textsc{FE}2 multiscale approach for modelling the elastoviscoplastic
  behaviour of long fibre sic/ti composite materials.
\newblock {\em Computer Methods in Applied Mechanics and Engineering},
  183(3-4):309--330, 2000.

\bibitem{Kouznetsova-etal2001}
V.~Kouznetsova, W.~A.~M. Brekelmans, and F.~P.~T. Baaijens.
\newblock An approach to micro-macro modeling of heterogeneous materials.
\newblock {\em Computational Mechanics}, 27(1):37--48, 2001.

\bibitem{E-Engquist-Li-Ren-VandenEijnden2007}
Weinan E, Bj{\"o}rn Engquist, Xiantao Li, Weiqing Ren, and Eric Vanden-Eijnden.
\newblock Heterogeneous multiscale methods : A review.
\newblock {\em Communications in Computational Physics}, 2(3):367--450, 2007.

\bibitem{Assyr2009}
Assyr Abdulle.
\newblock The finite element heterogeneous multiscale method: A computational
  strategy for multiscale \textsc{PDE}s.
\newblock {\em GAKUTO International Series Mathematical Sciences and
  Applications}, 31:133--181, 2009.

\bibitem{Assyr-etal2012}
Assyr Abdulle, Weinan E, Bj{\"o}rn Engquist, and Eric Vanden-Eijnden.
\newblock The heterogeneous multiscale method.
\newblock {\em Acta Numerica}, 21:1--87, 2012.

\bibitem{AbdulleHuber2016}
Assyr Abdulle and Martin~E. Huber.
\newblock Error estimates for finite element approximations of nonlinear
  monotone elliptic problems with application to numerical homogenization.
\newblock {\em Numerical Methods for Partial Differential Equations},
  32(3):955--969, 2016.

\bibitem{HenningOhlberger2015}
Patrick Henning and Mario Ohlberger.
\newblock Error control and adaptivity for heterogeneous multiscale
  approximations of nonlinear monotone problems.
\newblock {\em Discrete and Continuous Dynamical Systems - Series S},
  8(1):119--150, 2015.

\bibitem{Nejad-Wieners-2019}
Ramin~Shirazi Nejad and Christian Wieners.
\newblock Parallel inelastic heterogeneous multi-scale simulations.
\newblock In Stefan Diebels and Sergej Rjasanow, editors, {\em Multi-scale
  Simulation of Composite Materials}, volume~23 of {\em Mathematical
  Engineering}, pages 57--96. {Springer Berlin Heidelberg}, Berlin, Heidelberg,
  2019.

\bibitem{E-Ming-Zhang-2005}
Weinan E, Pingbing Ming, and Pingwen Zhang.
\newblock Analysis of the heterogeneous multiscale method for elliptic
  homogenization problems.
\newblock {\em Journal of the American Mathematical Society}, 18(01):121--157,
  2005.

\bibitem{Ohlberger2005}
Mario Ohlberger.
\newblock A posteriori error estimates for the heterogeneous multiscale finite
  element method for elliptic homogenization problems.
\newblock {\em Multiscale Modeling {\&} Simulation}, 4(1):88--114, 2005.

\bibitem{EidelFischer2018a}
Andreas Fischer and Bernhard Eidel.
\newblock Convergence and error analysis of \textsc{FE-HMM/FE2} for
  energetically consistent micro-coupling conditions in linear elastic solids.
\newblock {\em European Journal of Mechanics - A/Solids}, 77, 2019.

\bibitem{Moulinec-Suquet-1994}
H.~Moulinec and P.~Suquet.
\newblock A fast numerical method for computing the linear and nonlinear
  mechanical properties of composites.
\newblock {\em Comptes rendus de l'Acad\'{e}mie des Sciences Paris II},
  318:1417--14230, 1994.

\bibitem{Moulinec-Suquet-1998}
H.~Moulinec and P.~Suquet.
\newblock A numerical method for computing the overall response of nonlinear
  composites with complex microstructure.
\newblock {\em Computer Methods in Applied Mechanics and Engineering},
  157(1-2):69--94, 1998.

\bibitem{Zeller-Dederichs}
R.~Zeller and P.~H. Dederichs.
\newblock Elastic constants of polycrystals.
\newblock {\em Physica Status Solidi (b)}, 55(2):831--842, 1973.

\bibitem{Kroener-BOOK-1971}
Ekkehart Kr{\"o}ner.
\newblock {\em Statistical continuum mechanics: Course held ... October 1971},
  volume no. 92 of {\em Courses and lectures / International Centre for
  Mechanical Sciences}.
\newblock Springer, Wien, 1972.

\bibitem{Schneider-Ospald-Kabel-2015}
Matti Schneider, Felix Ospald, and Matthias Kabel.
\newblock Computational homogenization of elasticity on a staggered grid.
\newblock {\em International Journal for Numerical Methods in Engineering},
  105(9):693--720, 2016.

\bibitem{Schneider-Merkert-Kabel-2017}
Matti Schneider, Dennis Merkert, and Matthias Kabel.
\newblock \textsc{FFT}-based homogenization for microstructures discretized by
  linear hexahedral elements.
\newblock {\em International Journal for Numerical Methods in Engineering},
  109(10):1461--1489, 2017.

\bibitem{Yvonet-He-2007}
J.~Yvonnet and Q.-C. He.
\newblock The reduced model multiscale method \textsc{(R3M)} for the non-linear
  homogenization of hyperelastic media at finite strains.
\newblock {\em Journal of Computational Physics}, 223(1):341--368, 2007.

\bibitem{Radermacher-Reese}
Annika Radermacher and Stefanie Reese.
\newblock \textsc{POD}-based model reduction with empirical interpolation
  applied to nonlinear elasticity.
\newblock {\em International Journal for Numerical Methods in Engineering},
  107(6):477--495, 2016.

\bibitem{Soldner-etal2017}
Dominic Soldner, Benjamin Brands, Reza Zabihyan, Paul Steinmann, and Julia
  Mergheim.
\newblock A numerical study of different projection-based model reduction
  techniques applied to computational homogenisation.
\newblock {\em Computational Mechanics}, 60(4):613--625, 2017.

\bibitem{SchroederBalzaniBrands2010}
J{\"o}rg Schr{\"o}der, Daniel Balzani, and Dominik Brands.
\newblock Approximation of random microstructures by periodic statistically
  similar representative volume elements based on lineal-path functions.
\newblock {\em Archive of Applied Mechanics}, 81(7):975--997, 2011.

\bibitem{RheinbachKlawonn2006}
Axel Klawonn and Oliver Rheinbach.
\newblock Robust feti-dp methods for heterogeneous three dimensional elasticity
  problems.
\newblock {\em Computer Methods in Applied Mechanics and Engineering},
  196(8):1400--1414, 2007.

\bibitem{RheinbachKlawonn2010}
A.~Klawonn and O.~Rheinbach.
\newblock Highly scalable parallel domain decomposition methods with an
  application to biomechanics.
\newblock {\em Journal of Applied Mathematics and Mechanics}, 90(1):5--32,
  2010.

\bibitem{Holzer.2011}
Lorenz Holzer and Marco Cantoni.
\newblock Review of \textsc{FIB}-tomography.
\newblock In Ivo Utke, Stanislav Moshkalev, and Phillip Russell, editors, {\em
  Nanofabrication Using Focused Ion and Electron Beams}, volume 559201222 of
  {\em Oxford series on nanomanufacturing}, pages 410--435. {Oxford University
  Press}, Oxford, 2011.

\bibitem{Saenger-etal-2011}
Erik~H. Saenger, Frieder Enzmann, Youngseuk Keehm, and Holger Steeb.
\newblock Digital rock physics: Effect of fluid viscosity on effective elastic
  properties.
\newblock {\em Journal of Applied Geophysics}, 74(4):236--241, 2011.

\bibitem{Andrae-etal-2013}
Heiko Andr{\"a}, Nicolas Combaret, Jack Dvorkin, Erik Glatt, Junehee Han,
  Matthias Kabel, Youngseuk Keehm, Fabian Krzikalla, Minhui Lee, Claudio
  Madonna, Mike Marsh, Tapan Mukerji, Erik~H. Saenger, Ratnanabha Sain, Nishank
  Saxena, Sarah Ricker, Andreas Wiegmann, and Xin Zhan.
\newblock Digital rock physics benchmarks---\textsc{P}art \textsc{I}: Imaging
  and segmentation.
\newblock {\em Computers {\&} Geosciences}, 50:25--32, 2013.

\bibitem{Legrain-etal2011}
G.~Legrain, P.~Cartraud, I.~Perreard, and N.~Mo{\"e}s.
\newblock An \textsc{X-FEM} and level set computational approach for
  image-based modelling: Application to homogenization.
\newblock {\em International Journal for Numerical Methods in Engineering},
  86(7):915--934, 2011.

\bibitem{Lian-etal2013}
W.~D. Lian, G.~Legrain, and P.~Cartraud.
\newblock Image-based computational homogenization and localization: comparison
  between \textsc{X-FEM}/levelset and voxel-based approaches.
\newblock {\em Computational Mechanics}, 51(3):279--293, 2013.

\bibitem{FischerEidel2019}
Andreas Fischer and Bernhard Eidel.
\newblock Error analysis for quadtree-type mesh-coarsening algorithms adapted
  to pixelized heterogeneous microstructures, 2019.

\bibitem{Miehe-Koch-2002}
C.~Miehe and A.~Koch.
\newblock Computational micro-to-macro transitions of discretized
  microstructures undergoing small strains.
\newblock {\em Archive of Applied Mechanics (Ingenieur Archiv)},
  72(4-5):300--317, 2002.

\bibitem{Miehe-2003}
Christian Miehe.
\newblock Computational micro-to-macro transitions for discretized
  micro-structures of heterogeneous materials at finite strains based on the
  minimization of averaged incremental energy.
\newblock {\em Computer Methods in Applied Mechanics and Engineering},
  192(5-6):559--591, 2003.

\bibitem{Feyel2003}
Fr{\'e}d{\'e}ric Feyel.
\newblock A multilevel finite element method \textsc{(FE2)} to describe the
  response of highly non-linear structures using generalized continua.
\newblock {\em Computer Methods in Applied Mechanics and Engineering},
  192(28-30):3233--3244, 2003.

\bibitem{Kouznetsova-etal2004}
V.~G. Kouznetsova, M.G.D. Geers, and W.A.M. Brekelmans.
\newblock Multi-scale second-order computational homogenization of multi-phase
  materials: a nested finite element solution strategy.
\newblock {\em Computer Methods in Applied Mechanics and Engineering},
  193(48-51):5525--5550, 2004.

\bibitem{OezdemirBrekelmansGeers2008}
I.~{\"O}zdemir, W.A.M. Brekelmans, and M.G.D. Geers.
\newblock \textsc{FE}$^2$ computational homogenization for the
  thermo-mechanical analysis of heterogeneous solids.
\newblock {\em Computer Methods in Applied Mechanics and Engineering},
  198(3-4):602--613, 2008.

\bibitem{TemizerWriggers2008}
{\.{I}}.~Temizer and P.~Wriggers.
\newblock On the computation of the macroscopic tangent for multiscale
  volumetric homogenization problems.
\newblock {\em Computer Methods in Applied Mechanics and Engineering},
  198(3-4):495--510, 2008.

\bibitem{Schroeder2014}
J{\"o}rg Schr{\"o}der.
\newblock A numerical two-scale homogenization scheme: the \textsc{FE}2-method.
\newblock In Friedrich Pfeiffer, Franz~G. Rammerstorfer, Elisabeth Guazzelli,
  Bernhard Schrefler, Paolo Serafini, J{\"o}rg Schr{\"o}der, and Klaus Hackl,
  editors, {\em Plasticity and Beyond}, volume 550 of {\em CISM International
  Centre for Mechanical Sciences}, pages 1--64. {Springer Vienna}, Vienna,
  2014.

\bibitem{Saeb-Steinmann-Javili-2016}
Saba Saeb, Paul Steinmann, and Ali Javili.
\newblock Aspects of computational homogenization at finite deformations: A
  unifying review from \textsc{R}euss' to \textsc{V}oigt's bound.
\newblock {\em Applied Mechanics Reviews}, 68(5):050801, 2016.

\bibitem{KlawonnKohelerLanserRheinbach2019}
Axel Klawonn, Stephan K{\"o}hler, Martin Lanser, and Oliver Rheinbach.
\newblock Computational homogenization with million-way parallelism using
  domain decomposition methods.
\newblock {\em Computational Mechanics}, 21(1):1, 2019.

\bibitem{Eidel-etal-2018}
Bernhard Eidel, Andreas Fischer, and Ajinkya Gote.
\newblock A nonlinear \textsc{FE-HMM} formulation along with a novel
  algorithmic structure for finite deformation elasticity.
\newblock {\em PAMM}, 18(1):e201800457, 2018.

\bibitem{Bensoussan-Lions-Papanicolau-BOOK-1976}
G.~Papanicolau, A.~Bensoussan, and J.~L. Lions.
\newblock {\em Asymptotic Analysis for Periodic Structures}.
\newblock Elsevier, Burlington, 1978.

\bibitem{Sanchez-Palencia-BOOK-1980}
Enrique Sanchez-Palencia.
\newblock {\em Non-homogenous media and vibration theory}.
\newblock Springer-Verlag, Berlin, 1980.

\bibitem{GuedesKikuchi1990}
Jos{\'e}Miranda Guedes and Noboru Kikuchi.
\newblock Preprocessing and postprocessing for materials based on the
  homogenization method with adaptive finite element methods.
\newblock {\em Computer Methods in Applied Mechanics and Engineering},
  83(2):143--198, 1990.

\bibitem{Allaire1992}
Gr{\'e}goire Allaire.
\newblock Homogenization and two-scale convergence.
\newblock {\em SIAM Journal on Mathematical Analysis}, 23(6):1482--1518, 1992.

\bibitem{Cioranescu-Donato-BOOK-1999}
D.~Cioranescu and Patrizia Donato.
\newblock {\em An introduction to homogenization}, volume~17 of {\em Oxford
  lecture series in mathematics and its applications}.
\newblock {Oxford University Press}, Oxford and New York, 1999.

\bibitem{Hill1963}
R.~Hill.
\newblock Elastic properties of reinforced solids: Some theoretical principles.
\newblock {\em Journal of the Mechanics and Physics of Solids}, 11(5):357--372,
  1963.

\bibitem{Ostoja-Starzewski2006}
Martin Ostoja-Starzewski.
\newblock Material spatial randomness: From statistical to representative
  volume element.
\newblock {\em Probabilistic Engineering Mechanics}, 21(2):112--132, 2006.

\bibitem{Sab-1992}
K.~Sab.
\newblock On the homogenization and the simulation of random materials.
\newblock {\em European Journal of Mechanics A-Solids}, 11(5):585--607, 1992.

\bibitem{Drugan-Willis-1996}
W.~J. Drugan and J.~R. Willis.
\newblock A micromechanics-based nonlocal constitutive equation and estimates
  of representative volume element size for elastic composites.
\newblock {\em Journal of the Mechanics and Physics of Solids}, 44(4):497--524,
  1996.

\bibitem{Mandel-BOOK-1971}
Jean Mandel.
\newblock {\em Plasticit\'{e} classique et viscoplasticit\'{e}}, volume No 97
  of {\em Courses and Lectures, International centre for mechanical sciences}.
\newblock Springer, Vienna, 1972.

\bibitem{Hill1972}
R.~Hill.
\newblock On constitutive macro-variables for heterogeneous solids at finite
  strain.
\newblock {\em Proceedings of the Royal Society A: Mathematical, Physical and
  Engineering Sciences}, 326(1565):131--147, 1972.

\bibitem{E-Engquist-2003}
Weinan E and Bjorn Engquist.
\newblock The heterogeneous multiscale methods.
\newblock {\em Communications in Mathematical Sciences}, 1(1):87--132, 2003.

\bibitem{SouzaNeto-Feijoo2008}
E.~A. de~{Souza Neto} and R.~A. Feij{\'o}o.
\newblock On the equivalence between spatial and material volume averaging of
  stress in large strain multi-scale solid constitutive models.
\newblock {\em Mechanics of Materials}, 40(10):803--811, 2008.

\bibitem{GrytzMeschke2008}
R.~Grytz and G.~Meschke.
\newblock Consistent micro-macro transitions at large objective strains in
  curvilinear convective coordinates.
\newblock {\em International Journal for Numerical Methods in Engineering},
  73(6):805--824, 2008.

\bibitem{vDijk-2016}
N.~P. {van Dijk}.
\newblock Formulation and implementation of stress-driven and/or strain-driven
  computational homogenization for finite strain.
\newblock {\em International Journal for Numerical Methods in Engineering},
  107(12):1009--1028, 2016.

\bibitem{EidelFischer2016}
Bernhard Eidel and Andreas Fischer.
\newblock The heterogeneous multiscale finite element method \textsc{FE-HMM}
  for the homogenization of linear elastic solids.
\newblock {\em PAMM}, 16(1):521--522, 2016.

\bibitem{JeckerAbdulle2016}
Orane Jecker and Assyr Abdulle.
\newblock Numerical experiments for multiscale problems in linear elasticity.
\newblock In B{\"u}lent Karas{\"o}zen, Murat Manguo{\u{g}}lu, M{\"u}nevver
  Tezer-Sezgin, Serdar G{\"o}ktepe, and {\"O}m{\"u}r U{\u{g}}ur, editors, {\em
  Numerical Mathematics and Advanced Applications ENUMATH 2015}, volume 112 of
  {\em Lecture Notes in Computational Science and Engineering}, pages 123--131.
  {Springer International Publishing}, Cham, 2016.

\bibitem{Fujita-etal-2012}
Takeshi Fujita, Pengfei Guan, Keith McKenna, Xingyou Lang, Akihiko Hirata, Ling
  Zhang, Tomoharu Tokunaga, Shigeo Arai, Yuta Yamamoto, Nobuo Tanaka, Yoshifumi
  Ishikawa, Naoki Asao, Yoshinori Yamamoto, Jonah Erlebacher, and Mingwei Chen.
\newblock Atomic origins of the high catalytic activity of nanoporous gold.
\newblock {\em Nature materials}, 11(9):775--780, 2012.

\bibitem{Kramer-etal2004}
Dominik Kramer, Raghavan~Nadar Viswanath, and J{\"o}rg Weissm{\"u}ller.
\newblock Surface-stress induced macroscopic bending of nanoporous gold
  cantilevers.
\newblock {\em Nano Letters}, 4(5):793--796, 2004.

\bibitem{GriffithsBargmannReddy2017}
Emma Griffiths, Swantje Bargmann, and B.~D. Reddy.
\newblock Elastic behaviour at the nanoscale of innovative composites of
  nanoporous gold and polymer.
\newblock {\em Extreme Mechanics Letters}, 17:16--23, 2017.

\end{thebibliography}
